\documentclass{article}
\usepackage{amssymb,amsmath,amsfonts,authblk,amscd}
\usepackage{tikz}
\usetikzlibrary{calc,intersections,through,backgrounds}

\usepackage[french]{babel}
\usepackage[utf8]{inputenc}
\usepackage[T1]{fontenc}
\usepackage{proof}

\newtheorem{exemple}{Exemple}[section]
\newtheorem{definition}[exemple]{D\'efinition}
\newtheorem{theoreme}[exemple]{Th\'eor\`eme}
\newtheorem{lemme}[exemple]{Lemme}
\newtheorem{corollaire}[exemple]{Corollaire}

\newtheorem{remarque}[exemple]{Remarque}
\newtheorem{conjecture}[exemple]{Conjecture}

\newtheorem{proposition}[exemple]{Proposition}

\title{Quelques méthodes pour les mots sturmiens}

\author{Anna Frid}
\date{Written for the Spring School in Mathematical Computer Science
(École jeunes chercheurs en informatique mathématique),
Marseille, 4 - 8 March 2019}
\begin{document}
\maketitle 

\begin{abstract}
Les mots sturmiens sont une famille classique de mots infinis apériodiques qui pour tout $n$ possèdent seulement seulement $n + 1$ facteurs distincts de longueur n, étant la valeur minimale pour tout mot apériodique. Nous discuterons de leurs définitions équivalentes et plusieurs méthodes pour leur traitement qui, à ma connaissance, n'ont pas encore été décrites dans un livre. En particulier, nous considérerons la dualité de Berstel et Pocchiola ainsi que le lien entre les facteurs palindromes et le système de numération d'Ostrowski. Nous discuterons aussi la conjecture sur la longueur palindromique, qui a été démontrée sur les mots sturmiens par deux méthodes complètement différentes pour deux cas distincts.
\end{abstract}

\section{Notations et définitions équivalentes}

La famille des mots sturmiens est si classique, bien étudiée et toujours non triviale, si bien qu’il est très difficile d’écrire un
 nouveau texte d'introduction sur elle. D'une part, il existe plusieurs chapitres de monographies, importants et bien écrits, notamment \cite{lothaire}, \cite{fogg} et \cite{a_sh}. D'autre part, les preuves de nombreuses propriétés importantes restent compliquées et prennent de longues pages techniques. Je vais donc me concentrer sur des aspects non pris en compte dans les monographies précédentes, y compris la technique de dualité de Berstel et Pocchiola pour le calcul de la complexité totale et le lien entre les occurrences de palindromes et les systèmes de numération d'Ostrowski. Pour se concentrer sur cela, on doit laisser de nombreuses propositions avec des références aux papiers originaux au lieu des démonstrations. Pour le compenser, je vais essayer d'expliquer comment utiliser plusieurs techniques qui ne sont pas encore courantes.

Comme d'habitude, nous considérons des mots finis et infinis sur un alphabet fini; le seul cas que nous considérons est l’alphabet binaire, il est parfois commode de le désigner comme $\{0,1\}$ et parfois comme $\{a, b\} $. Un mot fini est noté $u = u[1] \cdots u[n]$, où $u[i]$ sont des symboles de l'alphabet; un mot infini est noté ${\bf u} = {\bf u}[0] {\bf u}[1] \cdots {\bf u}[n] \cdots$ ou bien ${\bf u} = {\bf u}[1] {\bf u}[2] \cdots {\bf u}[ n] \cdots$, ici ${\bf u}[i]$ sont toujours des symboles de l'alphabet. Le choix de zéro ou un comme premier indice est déterminé par la méthode considérée. (Notez que c’est une bonne réponse à la question «Quel est le premier nombre naturel ?»: zéro, tel qu’il est accepté dans la tradition française, ou un, comme il est habituel dans la littérature anglophone ? La bonne réponse dépendra du problème considéré.)

Un {\it facteur} d'un mot fini $u$ ou infini  ${\bf u}$ est tout simplement un mot fini de la forme 
$u[i] u[i + 1] \cdots u[j-1] u[j]$, où $j \geq i$; nous disons également que le mot vide $\epsilon$ ne contenant aucun symbole est un facteur de tout mot. Un facteur $u[i] u[i + 1] \cdots u[j-1] u[j]$ est noté aussi  $u[i..j]$, ou $u(i-1..j]$, ou 
$u[i..j + 1)$ (le symbole ${\bf u}$ est en gras si le mot est infini). L'ensemble des facteurs de ${\bf u}$ est désigné par Fac$({\bf u})$.

La longueur d'un mot fini $u = u[1] \cdots u[n]$ est définie comme $n$ et notée $|u|$. Le nombre d'occurrences d'une lettre $x$ dans $u$ est noté $|u|_x$. On dit qu'un mot infini ${\bf u}$ est {\it ultimement ($t-$)périodique} si 
${\bf u}[i] = {\bf u}[i + t]$ pour tout $i \geq N$ pour un certain $N$; un mot qui n'est pas ultimement périodique est appelé {\it apériodique}.

La {\it complexité} ({\it combinatoire}, {\it factor complexity} en anglais) $p_{\bf u}(n)$ d'un mot infini ${\bf u}$ est la fonction égale pour chaque $n$ au nombre de facteurs de ${\bf u}$ de longueur $n$.

\begin{exemple}\label{e:e1}
 {\rm
Le mot infini
\[{\bf u}=aababababababab \cdots\]
est ultimement 2-périodique. On a Fac$({\bf u})=\{a,b,aa,ab,ba,aab,aba,bab,\cdots\}$, et donc $p_{\bf u}(1)=2$ et $p_{\bf u}(n)=3$ pour tout $n>1$.

Le mot infini
\[{\bf w}=aabaabaaabaabaaabaabaabaaabaabaaabaabaaba\cdots\]
ne semble pas périodique, et on voit 
\[{\mbox Fac}({\bf w}){=}\{a,b,aa,ab,ba,aab,aba,baa,aaa,aaab,aaba,abaa,baaa,baab,\cdots\}.\] 
Dans la partie observée du mot, on a donc $p_{\bf w}(1)=2$, $p_{\bf w}(2)=3$, $p_{\bf w}(3)=4$ et $p_{\bf w}(4)=5$.
 }
\end{exemple}

Les résultats du reste de cette section sont discutés et rigoureusement démontrés dans \cite{lothaire}.

Un théorème classique et simple de Morse et Hedlund, daté de 1940, dit que la plus petite complexité d’un mot infini ${\bf u}$ qui n’est pas ultimement périodique est $p_{\bf u}(n) = n + 1$. Si cette égalité est vraie pour tout $n$, on dit que le mot est {\it sturmien}. La première question qui se pose est bien sûr l’existence des mots sturmiens: même si le mot ${\bf w}$ de l’exemple précédent ressemble à un mot sturmien, pouvons-nous toujours en construire la continuation sturmienne? La réponse est positive, et on sait bien construire des mots sturmiens.

Étant donné une {\it pente} $\sigma \in (0,1)$ et un {\it intercept} $\rho$, également sans perte de généralité appartenant à 
$[0,1)$, nous définissons le {\it mot mécanique inférieur}
${\bf s}_{\sigma, \rho}$ par les égalités
\[{\bf s}[n] = \lfloor (n + 1) \sigma + \rho \rfloor - \lfloor n \sigma + \rho \rfloor \]
pour tout $n \geq 0$.

Symétriquement, nous définissons le {\it mot mécanique supérieur} ${\bf s '}_{\sigma, \rho}$ par les égalités
\[{\bf s '}[n] = \lceil (n + 1) \sigma + \rho \rceil - \lceil n \sigma + \rho \rceil. \]

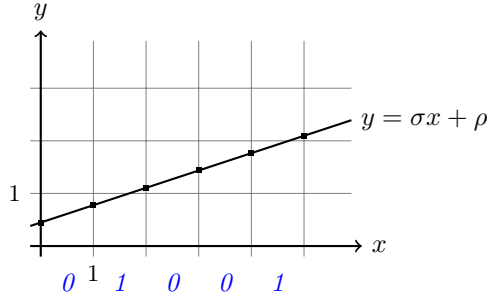
\begin{figure}
\centering
\begin{tikzpicture}[scale=.7]
  \begin{scope}
\draw[step=1cm,gray,very thin,name path=ggg] (-0.2,-0.2) grid (5.9,3.9);
 \draw[thick,->] (-0.2,0) -- (6.1,0) node[right] {$x$};
\draw[thick,->] (0,-0.2) -- (0,4.1) node[above] {$y$};
\draw (1,-0.2) node[below] {\small 1};
\draw (-0.2,1) node[left] {\small 1};
\draw [thick,name path=sigmarho](-0.2,0.381)--(5.9,2.392) node[right] {$y=\sigma x +\rho$};
\path [name path=Y0] (0,0)--(0,5);
\path [name intersections={of=sigmarho and Y0,by=E0}];
\node [fill=black,inner sep=1pt] at (E0) {};

\path [name path=Y1] (1,0)--(1,5);
\path [name intersections={of=sigmarho and Y1,by=E1}];
\node [fill=black,inner sep=1pt] at (E1) {};

\path [name path=Y2] (2,0)--(2,5);
\path [name intersections={of=sigmarho and Y2,by=E2}];
\node [fill=black,inner sep=1pt] at (E2) {};

\path [name path=Y3] (3,0)--(3,5);
\path [name intersections={of=sigmarho and Y3,by=E3}];
\node [fill=black,inner sep=1pt] at (E3) {};

\path [name path=Y4] (4,0)--(4,5);
\path [name intersections={of=sigmarho and Y4,by=E4}];
\node [fill=black,inner sep=1pt] at (E4) {};

\path [name path=Y5] (5,0)--(5,5);
\path [name intersections={of=sigmarho and Y5,by=E5}];
\node [fill=black,inner sep=1pt] at (E5) {};

\draw [color=blue] (0.5, -0.7) node {\it 0};
\draw [color=blue] (1.5, -0.7) node {\it 1};
\draw [color=blue] (2.5, -0.7) node {\it 0};
\draw [color=blue] (3.5, -0.7) node {\it 0};
\draw [color=blue] (4.5, -0.7) node {\it 1};

  \end{scope}
\end{tikzpicture}
\caption{Un mot sturmien comme un mot mécanique}\label{f:defst}
\end{figure}

La pente et l'intercept d'un mot mécanique sont bien la pente et l'interception de la droite $y = \sigma x + \rho $, voir la figure~\ref{f:defst}. Un zéro dans un mot mécanique correspond au cas où la droite coupe deux droites verticales  entières consécutives dans le même intervalle entier, et une unité signifie que la droite $y = \sigma x + \rho$ a croisé une ligne horizontale entière entre $x=n$ et $x=n+1$. Ainsi, un mot mécanique code et discrétise une ligne droite d'une pente quelconque.

\begin{figure} 
\centering
\begin{tikzpicture}
  \begin{scope}
\draw[step=1cm,gray,very thin,name path=ggg] (0.2,0.2) grid (1.8,1.8);
\draw [thick,name path=sigmarho](0.2,0.618)--(1.8,1.382) node[right] {};

\draw (1, 1) node {$\bullet$};
\draw (1, -0.2) node {\it m};
\draw  (-0.2, 1) node {\it n};

  \end{scope}

\end{tikzpicture}
\caption{Le cas de $\sigma m +\rho=n$}\label{f:mn}
\end{figure}
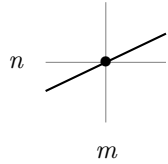

La différence entre les mots mécaniques inférieurs et supérieurs n'existe que lorsque la ligne $y = \sigma x + \rho$ passe par un point entier, c'est-à-dire lorsque $\sigma m + \rho \in \mathbb Z $ pour un $m \geq 0$ (voir Fig.~\ref{f:mn}). Dans ce cas, dans le mot mécanique inférieur, nous avons ${\bf s}[m-1] = 0$ et ${\bf s}[m] = 1$, et dans le mot mécanique supérieur, 
${\bf s'}[m-1] = 1$ et ${\bf s'}[m] = 0$.

L’ensemble des facteurs d’un mot mécanique ne dépend que du paramètre $\sigma$: ni $\rho$ ni le choix du mot inférieur ou supérieur ne sont importants. On peut désigner Fac$({\bf s}_{\sigma, \rho})$ = Fac$({\bf s '}_{\sigma, \rho})$ = Fac$(\sigma)$. On peut aussi, sant perte de généralité, considérer des mots mécaniques inférieurs et pas supérieurs dans tous les problèmes liés à l'ensemble des facteurs.

Si $\sigma$ est rationnel, les mots mécaniques respectifs sont tous périodiques. Si $\sigma$ est irrationnel, il n’est pas très difficile de démontrer que les mots mécaniques sont sturmiens, donnant ainsi une famille d’exemples de mots sturmiens. En fait, les mots sturmiens sur l'alphabet $\{0,1 \}$ sont {\it exactement} les mots mécaniques de pentes irrationnelles. On peut donc les définir comme les mots de complexité $n + 1$ ou comme les mots mécaniques apériodiques. Ces définitions sont équivalentes.

\begin{remarque}
 {
 L'équivalence n'est observée que pour des mots infinis unilatéraux. Par exemple, le mot bi-infini $\cdots 000010000 \cdots$ est de complexité $n + 1$, mais il n'est pas un mot mécanique apériodique.
}
\end{remarque}
 
La définition d'un mot mécanique inférieur peut également être récrite comme suit:
\begin{equation}\label{e:frpart}
{\bf s}[n]=\begin{cases}
              0, \mbox{~si~} \{n\sigma+\rho \}\leq 1-\sigma,\\
              1, \mbox{~sinon},
             \end{cases}
\end{equation}
où $\{x\}$ désigne la partie fractionnaire de $x$. Dans certains cas, nous préférerons cette définition. Notons également que puisque l'ensemble $\{n \sigma + \rho\}$ est dense sur $[0,1]$, alors
\[\mbox{Fac}(\sigma)=\{s_{\sigma,\rho}[0..n]|\rho \in [0,1[, n \in \mathbb N\}. \]

Une autre définition équivalente de mots sturmiens utilise leur propriété d'équilibre. On dit qu'un mot infini ${\bf w}$ est {\it équilibré} si pour tous facteurs $v$ et $w$ de $\bf w$ de même longueur et pour tout symbole $x$ de l'alphabet,
\[||v|_x-|w|_x|\leq 1.\]

Par exemple, dans le mot ${\bf w}$ de l'exemple \ref{e:e1}, il existe zéro ou une lettre $b$ dans un facteur de longueur $3$ et un ou deux $b$s dans un facteur de longueur $4$. S'il n'y a que deux valeurs possibles pour le nombre de $b$s (et donc de $a$s) dans les facteurs de chaque longueur $n$, alors ce mot infini est équilibré.

Un mot est sturmien si et seulement si il est apériodique et équilibré.

Ce sont les trois principales définitions équivalentes des mots sturmiens: ce sont les mots de complexité $n + 1$, les mots mécaniques apériodiques et les mots équilibrés apériodiques. La démonstration de cette équivalence prend plusieurs pages lourdes dans \cite{lothaire}. Il existe des dizaines d'autres définitions et caractérisations; dans la section \ref{s:directive}, nous en discuterons une, liée aux fractions continues.

Avant la fin de cette section, nous donnons encore quelques définitions classiques. Un {\it langage} est un ensemble de mots finis. Une langage est appelé {\it factoriel} s'il contient tous les facteurs de ses éléments. Evidemment, l'ensemble des facteurs d'un mot infini, ainsi que l'union de plusieurs langages factoriels, sont factoriels. Un langage factoriel est {\it prolongeable} si pour tout élément $u$, il existe des lettres $x$ et $y$ telles que $xuy$ est également élément du langage. Le langage des facteurs de tout mot sturmien et le langage des facteurs de tous les mots sturmiens sont factoriels et prolongeables.

On dit qu'un élément $u$ d'un langage factoriel $F$ sur un alphabet $\Sigma$ est {\it spécial à gauche} si $xu$ et $yu$ appartiennent à $F$ pour deux lettres différentes $x$ et $y$. Symétriquement, on dit que $u$ est {\it spécial à droite} si 
$ux$ et $uy$ appartiennent à $F$ pour deux lettres différentes $x$ et $y$. Comme le langage des facteurs d'un mot sturmien est prolongeable et contient exactement $n + 1$ facteurs de chaque longueur, il est évident de voir qu'il contient exactement un facteur spécial à gauche et un facteur spécial à droite de chaque longueur.

Un élément d'un langage factoriel spécial à la fois  à gauche et à droite s'appelle {\it bispécial}. Les mots bispéciaux ont un intérêt particulier et nous les discuterons dans la section \ref{s:directive}.

\section{Nombre total de facteurs sturmiens}
Comme nous en avons discuté, un mot sturmien contient exactement $n + 1$ facteurs de longueur $n$, et tous les mots sturmiens de même pente contiennent le même ensemble de facteurs. Mais quel est le nombre de ces facteurs si nous considérons tous les mots sturmiens ensemble, en d’autres termes, quelle est la {\it complexité totale} de l’ensemble des mots sturmiens? Comme aucun mot sturmien ne contient le facteur $0011$, ce n'est pas $2^n$, alors quel est l'ordre de croissance de cette fonction?

La première réponse précise à cette question a été obtenue par Lipatov \cite{lipatov} en 1982. Lipatov a publié son article en russe sans utiliser le terme «mot sturmien»: il a parlé de «collections binaires» et de «classes d'uniformité». En termes modernes, il a utilisé les propriétés équilibrées des mots sturmiens. Son article n’a pas été remarqué dans la communauté occidentale (ne publiez jamais de nouveaux résultats dans une langue autre que l’anglais!). Le résultat a donc été redécouvert en 1991 par Mignosi \cite{mignosi}. Je donne ici la preuve géométrique de Berstel et Pocchiola, datée de 1993 \cite {bp}.

\begin{theoreme}\label{t:totalSt}
Le nombre total $p_s(n)$ de facteurs sturmiens de longueur $n$ est
\[p_s(n)= 1+\sum_{q=1}^{n} \varphi(q) (n+1-q) = \frac{n^3}{\pi^2}+O(n^2\log n),\]
où $\varphi(q)$ est l'indicatrice d'Euler, égale au nombre d'entiers compris entre 1 et $q$ (inclus) et premiers avec $q$.
\end{theoreme}
\noindent {\sc Démonstration.}
Considérons l'ensemble $S_n$ de tous les facteurs de mots sturmiens de longueur $n$. Chacun de ses éléments est un préfixe 
${\bf s}_{\sigma,\rho}[0..n)$ d'un mot sturmien infini ${\bf s}_{\sigma, \rho}$, ou, plus précisément, d'un continuum de mots sturmiens correspondant à un ensemble de paires $(\sigma, \rho) \subset [0,1)\times [0,1)$. Nous pouvons dessiner dans le carré $[0,1]\times [0,1]$ l'ensemble des points $(\sigma, \rho)$ qui donnent un préfixe $u \in S_n$: même si $\sigma$ est rationnel, c’est-à-dire si le mot infini ${\bf s}_{\sigma, \rho}$ est périodique et non sturmien, alors il a le même préfixe que certains mots sturmiens avec des paramètres proches.

Un {\it arrangement} ${\mathcal A}_n$ d'ordre $n$ est le carré $[0,1]\times [0,1]$ divisé en zones correspondant à différents préfixes ${\bf s}_{\sigma, \rho} [0..n)$. Étudions à quoi ressemblent les arrangements.

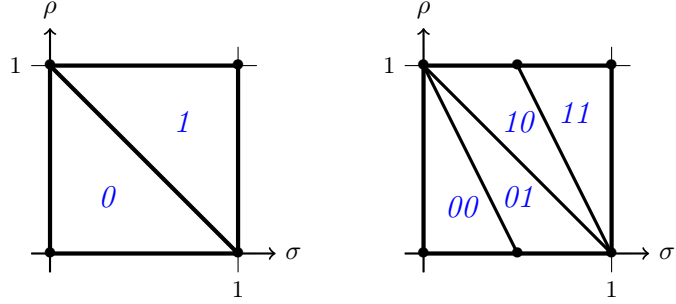
\begin{figure}
\centering
\begin{minipage}{0.4\textwidth}
\begin{tikzpicture}[scale=2.5]
  \begin{scope}
\draw[step=1cm, name path=ggg] (-0.1,-0.1) grid (1.1,1.1);
\draw[thick,->] (-0.1,0) -- (1.2,0) node[right] {$\sigma$};
\draw[thick,->] (0,-0.1) -- (0,1.2) node[above] {$\rho$};
\draw[ultra thick] (0,0) -- (0,1) -- (1,1) -- (1,0) -- cycle;
\draw[ultra thick] (0,1) -- (1,0);
\draw (1,-0.1) node[below] {\small 1};
\draw (-0.1,1) node[left] {\small 1};

\draw (1, 1) node {$\bullet$};
\draw (1, 0) node {$\bullet$};
\draw (0, 1) node {$\bullet$};
\draw (0, 0) node {$\bullet$};
\draw [color=blue] (0.3, 0.3) node {\it \large 0};
\draw [color=blue] (0.7, 0.7) node {\it \large 1};

  \end{scope}

\end{tikzpicture}
\end{minipage}
\begin{minipage}{0.4\textwidth}
 \begin{tikzpicture}[scale=2.5]
  \begin{scope}
\draw[step=1cm, name path=ggg] (-0.1,-0.1) grid (1.1,1.1);
\draw[thick,->] (-0.1,0) -- (1.2,0) node[right] {$\sigma$};
\draw[thick,->] (0,-0.1) -- (0,1.2) node[above] {$\rho$};
\draw[ultra thick] (0,0) -- (0,1) -- (1,1) -- (1,0) -- cycle;
\draw[very thick] (0,1) -- (1,0);
\draw[very thick] (0,1) -- (0.5,0);
\draw[very thick] (0.5,1) -- (1,0);
\draw (1,-0.1) node[below] {\small 1};
\draw (-0.1,1) node[left] {\small 1};

\draw (1, 1) node {$\bullet$};
\draw (1, 0) node {$\bullet$};
\draw (0, 1) node {$\bullet$};
\draw (0, 0) node {$\bullet$};
\draw (0.5, 0) node {$\bullet$};
\draw (0.5, 1) node {$\bullet$};
\draw [color=blue] (0.2, 0.25) node {\it \large 00};
\draw [color=blue] (0.5, 0.3) node {\it \large 01};
\draw [color=blue] (0.5, 0.7) node {\it \large 10};
\draw [color=blue] (0.8, 0.75) node {\it \large 11};

  \end{scope}

\end{tikzpicture}
\end{minipage}

\caption{Les arrangements d'ordre 1 et 2}\label{f:A1}
\end{figure}

La condition ${\bf s}[0] = \lfloor \sigma + \rho \rfloor - \lfloor \rho \rfloor = 1$ peut être immédiatement simplifiée comme $\sigma + \rho \geq 1$ puisque $ 0 \leq \rho <1$. Ainsi, la partie gauche de la figure~\ref{f:A1} montre l'arrangement 
${\mathcal A}_1$ avec la ligne de séparation $\sigma + \rho = 1$.

La condition ${\bf s}[1]=1$ veut dire exactement que $\lfloor 2\sigma+\rho \rfloor > \lfloor \sigma+\rho \rfloor$, c'est-à-dire, $\lfloor 2\sigma+\rho \rfloor = 1$ si $\lfloor \sigma+\rho \rfloor=0$ ou bien $\lfloor 2\sigma+\rho \rfloor = 2$ si $\lfloor \sigma+\rho \rfloor=1$. Pour obtenir ${\mathcal A}_2$ à partir de ${\mathcal A}_1$ nous devons donc ajouter les droites $2\sigma+\rho = 1$ et $2\sigma+\rho = 2$: ainsi, un mot sturmien de longueur $2$ est défini par la position du point   $(\sigma,\rho)$ par rapport à ces deux nouvelles droites ainsi que la droite précédente $\sigma+\rho=1$ (voir la partie droite de Fig.~\ref{f:A1}). 

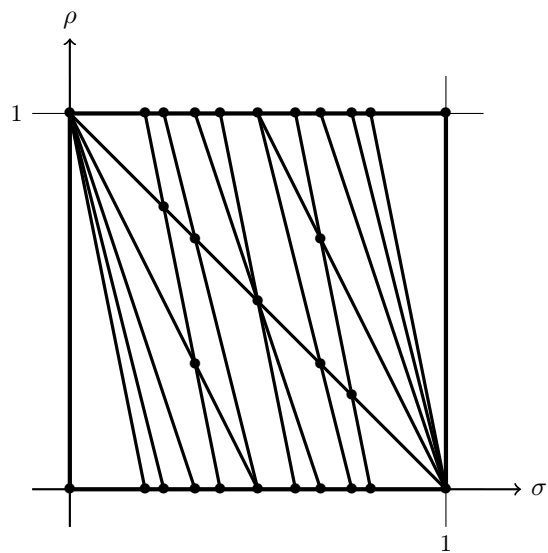
\begin{figure}
 \centering
\begin{tikzpicture}[scale=5]
  \begin{scope}
\draw[step=1cm, name path=ggg] (-0.1,-0.1) grid (1.1,1.1);
\draw[thick,->] (-0.1,0) -- (1.2,0) node[right] {$\sigma$};
\draw[thick,->] (0,-0.1) -- (0,1.2) node[above] {$\rho$};
\draw[ultra thick] (0,0) -- (0,1) -- (1,1) -- (1,0) -- cycle;
\draw[very thick] (0,1) -- (1,0);
\draw[very thick] (0,1) -- (0.5,0);
\draw[very thick] (0.5,1) -- (1,0);

\draw[very thick] (0,1) -- (0.333,0);
\draw[very thick] (0.333,1) -- (0.667,0);
\draw[very thick] (0.667,1) -- (1,0);

\draw[very thick] (0,1) -- (0.25,0);
\draw[very thick] (0.25,1) -- (0.5,0);
\draw[very thick] (0.5,1) -- (0.75,0);
\draw[very thick] (0.75,1) -- (1,0);

\draw[very thick] (0,1) -- (0.2,0);
\draw[very thick] (0.2,1) -- (0.4,0);
\draw[very thick] (0.4,1) -- (0.6,0);
\draw[very thick] (0.6,1) -- (0.8,0);
\draw[very thick] (0.8,1) -- (1,0);

\draw (1,-0.1) node[below] {\small 1};
\draw (-0.1,1) node[left] {\small 1};

\draw (1, 1) node {$\bullet$};
\draw (1, 0) node {$\bullet$};
\draw (0, 1) node {$\bullet$};
\draw (0, 0) node {$\bullet$};
\draw (0.5, 0) node {$\bullet$};
\draw (0.5, 1) node {$\bullet$};

\draw (0.333,1) node {$\bullet$};
\draw (0.667,1) node {$\bullet$};
\draw (0.333,0) node {$\bullet$};
\draw (0.667,0) node {$\bullet$};

\draw (0.5,0.5) node {$\bullet$};

\draw (0.25,1) node {$\bullet$};
\draw (0.75,1) node {$\bullet$};
\draw (0.25,0) node {$\bullet$};
\draw (0.75,0) node {$\bullet$};

\draw (0.66667,0.33333) node {$\bullet$};
\draw (0.33333,0.66667) node {$\bullet$};

\draw (0.2,1) node {$\bullet$};
\draw (0.4,1) node {$\bullet$};
\draw (0.6,1) node {$\bullet$};
\draw (0.8,1) node {$\bullet$};

\draw (0.2,0) node {$\bullet$};
\draw (0.4,0) node {$\bullet$};
\draw (0.6,0) node {$\bullet$};
\draw (0.8,0) node {$\bullet$};

\draw (0.66667,0.66667) node {$\bullet$};
\draw (0.33333,0.33333) node {$\bullet$};
\draw (0.25,0.75) node {$\bullet$};
\draw (0.75,0.25) node {$\bullet$};


  \end{scope}

\end{tikzpicture}
\caption{L'arrangement d'ordre 5}\label{f:A5}
\end{figure}

En continuant cet argument, nous voyons que l’arrangement ${\mathcal A}_n$ contient toutes les lignes de la forme $k\sigma+\rho=l$, où $k=1,\ldots,n$ et $l=1,\ldots,k$. Notons l'ensemble de leurs paramètres par $P_n$ et ajoutons à cet ensemble le point $(0,0)$:
\[P_n=\{(0,0)\}\bigcup_{1\leq l \leq k \leq n} \{(k,l)\}.\]
Un préfix sturmien ${\bf s}_{\sigma,\rho}[0..n)$ est défini par la position du point  $(\sigma,\rho)$ par rapport aux lignes droites avec des paramètres de l'ensemble $P_n$, c'est-à-dire, par la face de l'arrangement ${\mathcal A}_n$ qui contient le  point $(\sigma,\rho)$. Voici  l'arrangement ${\mathcal A}_5$ (voir Fig.~\ref{f:A5}).

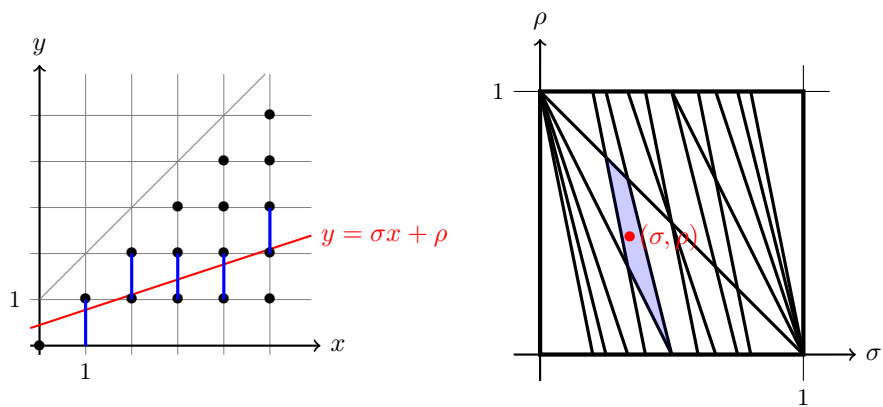
\begin{figure}
 \begin{minipage}{0.52\textwidth}
   \begin{tikzpicture}[scale=.6125]
     \begin{scope}
       \draw[step=1cm,gray,very thin,name path=ggg] (-0.2,-0.2) grid (5.9,5.9);
       \draw[thick,->] (-0.2,0) -- (6.1,0) node[right] {$x$};
       \draw[thick,->] (0,-0.2) -- (0,6.1) node[above] {$y$};
       \draw (1,-0.2) node[below] {\small 1};
       \draw (-0.2,1) node[left] {\small 1};
       \draw [gray, thin] (0,1) -- (4.9,5.9);
       \draw [thick,color=red,name path=sigmarho](-0.2,0.381)--(5.9,2.392) node[right] {$y=\sigma x +\rho$};

       \draw  (0,0) node {$\bullet$};
       \draw  (1,1) node {$\bullet$};
       \draw  (2,1) node {$\bullet$};
       \draw  (2,2) node {$\bullet$};
       \draw  (3,1) node {$\bullet$};
       \draw  (3,2) node {$\bullet$};
       \draw  (3,3) node {$\bullet$};
       \draw  (4,1) node {$\bullet$};
       \draw  (4,2) node {$\bullet$};
       \draw  (4,3) node {$\bullet$};
       \draw  (4,4) node {$\bullet$};
       \draw  (5,1) node {$\bullet$};
       \draw  (5,2) node {$\bullet$};
       \draw  (5,3) node {$\bullet$};
       \draw  (5,4) node {$\bullet$};
       \draw  (5,5) node {$\bullet$};

       \draw [very thick, color=blue] (1,0)--(1,1);
       \draw [very thick, color=blue] (2,1)--(2,2);
       \draw [very thick, color=blue] (3,1)--(3,2);
       \draw [very thick, color=blue] (4,1)--(4,2);
       \draw [very thick, color=blue] (5,2)--(5,3);

     \end{scope}
   \end{tikzpicture}
 \end{minipage}
 \hfill
 \begin{minipage}{0.47\textwidth}
   \begin{tikzpicture}[scale=3.5]
     \begin{scope}
       \draw[step=1cm, name path=ggg] (-0.1,-0.1) grid (1.1,1.1);
       \draw[thick,->] (-0.1,0) -- (1.2,0) node[right] {$\sigma$};
       \draw[thick,->] (0,-0.1) -- (0,1.2) node[above] {$\rho$};
       \draw[ultra thick] (0,0) -- (0,1) -- (1,1) -- (1,0) -- cycle;
       \draw[very thick] (0,1) -- (1,0);
       \draw[very thick] (0,1) -- (0.5,0);
       \draw[very thick] (0.5,1) -- (1,0);

       \draw[very thick] (0,1) -- (0.333,0);
       \draw[very thick] (0.333,1) -- (0.667,0);
       \draw[very thick] (0.667,1) -- (1,0);

       \draw[very thick] (0,1) -- (0.25,0);
       \draw[very thick] (0.25,1) -- (0.5,0);
       \draw[very thick] (0.5,1) -- (0.75,0);
       \draw[very thick] (0.75,1) -- (1,0);

       \draw[very thick] (0,1) -- (0.2,0);
       \draw[very thick] (0.2,1) -- (0.4,0);
       \draw[very thick] (0.4,1) -- (0.6,0);
       \draw[very thick] (0.6,1) -- (0.8,0);
       \draw[very thick] (0.8,1) -- (1,0);

       \draw (1,-0.1) node[below] {\small 1};
       \draw (-0.1,1) node[left] {\small 1};

       \fill [color=blue, opacity=0.2] (0.5,0) -- (0.333,0.667) -- (0.25,0.75) -- (0.333,0.333) -- cycle;

       \draw [color=red] (0.34,0.4469) node {$\bullet$} node[right] {$(\sigma,\rho)$};
     \end{scope}
   \end{tikzpicture}  
 \end{minipage}
 \caption{Dualité entre droites et points}\label{f:duality}
\end{figure}

Si $(\sigma,\rho)$ se trouve sur l'une des lignes avec les paramètres de $P_n$, cela correspond au cas considéré sur la Fig.~\ref{f:mn}, et de toute façon, ${\bf s}_{\sigma,\rho}[0..n)$ est égal à un préfixe sturmien correspondant à la partie interne d'une des faces voisines. Il suffit donc de considérer la partie interne des faces. De plus, non seulement une ligne  $y=\sigma x + \rho$ correspond à un point $(\sigma,\rho)$ de l'arrangement, mais une droite $k \sigma + \rho = l$ sur l'arrangement correspond au point $(k,l)\in P_n$ de la grille entière initiale. En effet, la position d'un point $(\sigma,\rho)$ par rapport aux droites $k \sigma + \rho = l$ correspond exactement à la position de la droite $y=\sigma x + \rho$ par rapport aux points de $P_n$ (voir Fig.~\ref{f:duality}). Cela implique notamment que des faces différentes de l’arrangement correspondent à des préfixes sturmiens différents, car deux positions différentes de la ligne $y=\sigma x + \rho$ par rapport aux points entiers $(k,l)$ donnent des 
mots sturmiens différents. C'est pourquoi
Berstel et Pocchiola appellent
la transformation $y=\sigma x + \rho  \leftrightarrow (\sigma,\rho)$ la {\it transformation de dualité}.

Ainsi, le nombre $p_s(n)$ de mots sturmiens de longueur $n$ est égal au nombre de faces de l'arrangement ${\mathcal A}_n$. Interprétons un arrangement comme un graphe planaire et calculons le nombre de faces qu’il contient avec la formule d’Euler
\[v-e+f=1.\]

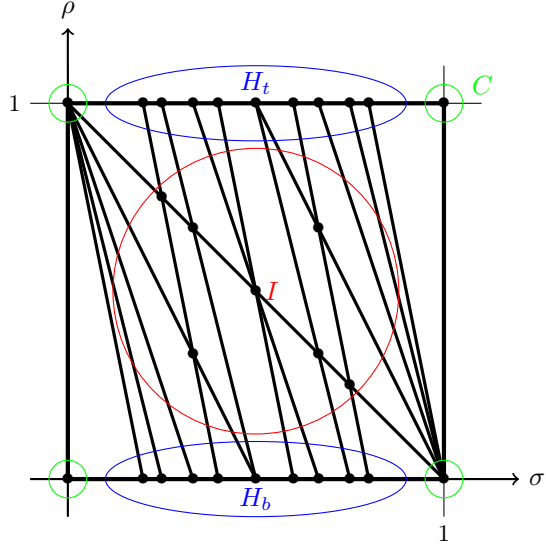
\begin{figure}
 \centering

\begin{tikzpicture}[scale=5]
  \begin{scope}
\draw[step=1cm, name path=ggg] (-0.1,-0.1) grid (1.1,1.1);
\draw[thick,->] (-0.1,0) -- (1.2,0) node[right] {$\sigma$};
\draw[thick,->] (0,-0.1) -- (0,1.2) node[above] {$\rho$};
\draw[ultra thick] (0,0) -- (0,1) -- (1,1) -- (1,0) -- cycle;
\draw[very thick] (0,1) -- (1,0);
\draw[very thick] (0,1) -- (0.5,0);
\draw[very thick] (0.5,1) -- (1,0);

\draw[very thick] (0,1) -- (0.333,0);
\draw[very thick] (0.333,1) -- (0.667,0);
\draw[very thick] (0.667,1) -- (1,0);

\draw[very thick] (0,1) -- (0.25,0);
\draw[very thick] (0.25,1) -- (0.5,0);
\draw[very thick] (0.5,1) -- (0.75,0);
\draw[very thick] (0.75,1) -- (1,0);

\draw[very thick] (0,1) -- (0.2,0);
\draw[very thick] (0.2,1) -- (0.4,0);
\draw[very thick] (0.4,1) -- (0.6,0);
\draw[very thick] (0.6,1) -- (0.8,0);
\draw[very thick] (0.8,1) -- (1,0);

\draw (1,-0.1) node[below] {\small 1};
\draw (-0.1,1) node[left] {\small 1};

\draw (1, 1) node {$\bullet$};
\draw (1, 0) node {$\bullet$};
\draw (0, 1) node {$\bullet$};
\draw (0, 0) node {$\bullet$};
\draw (0.5, 0) node {$\bullet$};
\draw (0.5, 1) node {$\bullet$};

\draw (0.333,1) node {$\bullet$};
\draw (0.667,1) node {$\bullet$};
\draw (0.333,0) node {$\bullet$};
\draw (0.667,0) node {$\bullet$};

\draw (0.5,0.5) node {$\bullet$};

\draw (0.25,1) node {$\bullet$};
\draw (0.75,1) node {$\bullet$};
\draw (0.25,0) node {$\bullet$};
\draw (0.75,0) node {$\bullet$};

\draw (0.66667,0.33333) node {$\bullet$};
\draw (0.33333,0.66667) node {$\bullet$};

\draw (0.2,1) node {$\bullet$};
\draw (0.4,1) node {$\bullet$};
\draw (0.6,1) node {$\bullet$};
\draw (0.8,1) node {$\bullet$};

\draw (0.2,0) node {$\bullet$};
\draw (0.4,0) node {$\bullet$};
\draw (0.6,0) node {$\bullet$};
\draw (0.8,0) node {$\bullet$};

\draw (0.66667,0.66667) node {$\bullet$};
\draw (0.33333,0.33333) node {$\bullet$};
\draw (0.25,0.75) node {$\bullet$};
\draw (0.75,0.25) node {$\bullet$};


\draw [color=green] (0,1) circle (0.05);
\draw [color=green] (1,1) circle (0.05);
\draw [color=green] (0,0) circle (0.05);
\draw [color=green] (1,0) circle (0.05);
\draw [color=green](1.05,1.05) node[right] {$C$};

\draw [color=blue] (0.5,0) ellipse (0.4 and 0.1) node[below] {$H_b$};
\draw [color=blue] (0.5,1) ellipse (0.4 and 0.1) node[above] {$H_t$};

\draw [color=red] (0.5,0.5) circle (0.38) node[right] {$I$};

  \end{scope}

\end{tikzpicture}

\caption{Quatre groupes de sommets de l'arrangement}\label{f:CHI}
\end{figure}

Ici $v$ est le nombre de sommets du graphe, $e$ est le nombre de ses arêtes, et le nombre de faces $f$ est le nombre $p_s(n)$ que nous recherchons. La somme $v-e + f$ est égale à $1$ puisque la face externe du graphe ne correspond pas à un mot sturmien et est donc exclue. Les sommets du graphe sont les quatre angles de l'arrangement plus toutes les intersections des lignes 
$k\sigma + \rho = l$, où $1 \leq l \leq k \leq n$, entre elles et avec les lignes horizontales $\rho = 0$ et $\rho = 1$: Voir la Fig.~\ref{f:CHI} où ces trois types de sommets sont marqués par $C$ pour les angles, $H_t$ et $H_b$ pour les intersections avec les deux lignes horizontales et $I$ pour les intersections internes entre deux lignes inclinées. Le nombre de lignes inclinées ou horizontales se croisant en un point $t$ est noté $c(t)$, donc par exemple pour $t = (0.5, 0.5)$ à la Fig.~\ref{f:CHI}, nous avons $c(t) = 3$.

Les arêtes du graphe sont les segments des droites entre deux points consécutifs. Comme d'habitude, le nombre d'arêtes $e$ d'un graphe est égal à la moitié de la somme des degrés de ses sommets:
\[p_s(n)=f=1+e-v= 1+\sum_{t\in V}\left(\frac{deg(t)}{2}-1\right),\]
où $V=C \cup H_t \cup H_b \cup I$.

Étudions les degrés des sommets de l'arrangement. Tout d'abord, il y a deux coins de degré 2 et deux coins de degré $n + 2$. Les ensembles $H_t$ et $H_b$ sont symétriques et le degré de tout sommet $t$ de l'un d'eux est $c(t) + 1$. Le degré d'un sommet $t$ de l'ensemble interne $I$ est $2c(t)$. En résumé, nous obtenons
\begin{align}
p_s(n)&=1+\sum_{t\in C}\left(\frac{deg(t)}{2}-1\right)+\sum_{t\in H_t\cup H_b}\left(\frac{deg(t)}{2}-1\right)+\sum_{t\in I}\left(\frac{deg(t)}{2}-1\right) \nonumber \\
&=1+n+2\sum_{t\in H_b }\left(\frac{deg(t)}{2}-1\right)+\sum_{t\in I}\left(\frac{deg(t)}{2}-1\right) \nonumber \\
&=1+n+\sum_{t\in H_b \cup I }\left(c(t)-1\right). \label{e1}
\end{align}

Pour calculer la somme ci-dessus, lions les sommets de $H_b \cup I$ aux lignes droites de la construction mécanique initiale. Chaque tel sommet $t$ est le point d'intersection de $c(t) \geq 2$ sur $n(n + 1)/2 + 1$ droites avec les paramètres de l'ensemble $P_n$: notons cet ensemble de lignes $L_n$. Sur l'image initiale, l'ensemble $L_n$ de droites correspond à l'ensemble $P_n$. Ainsi, un point $t = (\sigma, \rho)$ correspond à une droite $t^*: y = \sigma x + \rho$ qui passe par $c(t)$ points entiers à gauche de la Fig.~\ref{f:duality}.

\begin{figure}
 \centering
\begin{tikzpicture}[scale=.7]
  \begin{scope}
\draw[step=1cm,gray,very thin,name path=ggg] (-0.2,-0.2) grid (5.9,5.9);
 \draw[thick,->] (-0.2,0) -- (6.1,0) node[right] {$x$};
\draw[thick,->] (0,-0.2) -- (0,6.1) node[above] {$y$};
\draw (1,-0.2) node[below] {\small 1};
\draw (-0.2,1) node[left] {\small 1};
\draw [gray, thin] (0,1) -- (4.9,5.9);
\draw [very thick,color=blue, name path=sigmarho](-0.2,0.4)--(5.9,3.45); 

\draw  (0,0) node {$\bullet$};
\draw  (1,1) node {$\bullet$};
\draw  (2,1) node {$\bullet$};
\draw  (2,2) node {$\bullet$};
\draw  (3,1) node {$\bullet$};
\draw  (3,2) node {$\bullet$};
\draw  (3,3) node {$\bullet$};
\draw  (4,1) node {$\bullet$};
\draw  (4,2) node {$\bullet$};
\draw  (4,3) node {$\bullet$};
\draw  (4,4) node {$\bullet$};
\draw  (5,1) node {$\bullet$};
\draw  (5,2) node {$\bullet$};
\draw  (5,3) node {$\bullet$};
\draw  (5,4) node {$\bullet$};
\draw  (5,5) node {$\bullet$};


  \end{scope}
\end{tikzpicture}
\caption{Une droite $t^*: y=0.5x+0.5$}\label{f:0505}
\end{figure}
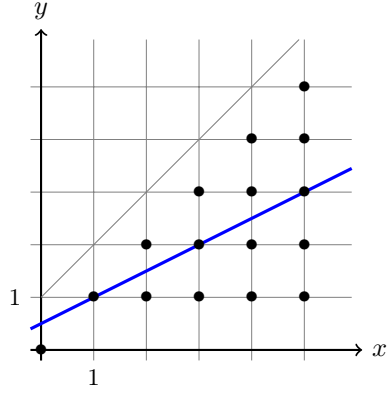

L'intersection $t$ des droites $k \sigma + \rho = l$ et $k'\sigma + \rho = l'$, $k> k'$, est le point
$t = \left (\frac{l-l'}{k-k'}, l-k \frac{l-l '}{k-k'} \right)$, dual à la droite 
$t^*: y = \frac{l-l'}{k-k'} x + l-k\frac{l-l'}{k-k'}$ de la construction initiale. Le paramètre $c(t)$ est égal au nombre de points de $P_n$ rencontrés par la droite $t^*$. Par exemple, le point $t = (0.5,0.5)$ correspond à la droite $y = 0,5 x + 0,5$ représentée à la Fig.~\ref{f:0505}.

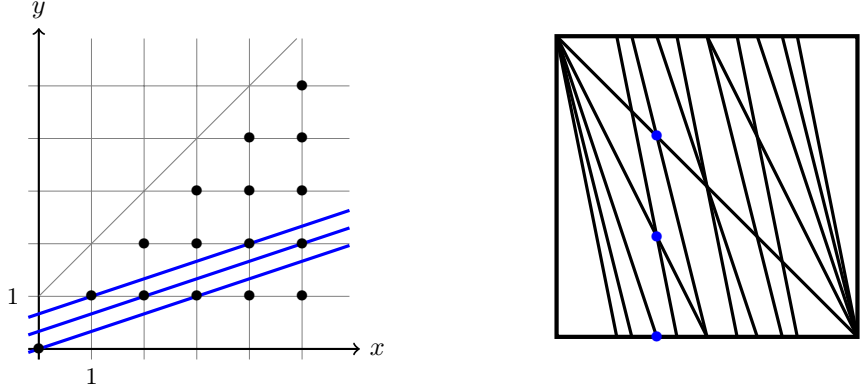
\begin{figure}
\centering
\begin{minipage}{0.6\textwidth}

 \begin{tikzpicture}[scale=.7]
  \begin{scope}
\draw[step=1cm,gray,very thin,name path=ggg] (-0.2,-0.2) grid (5.9,5.9);
 \draw[thick,->] (-0.2,0) -- (6.1,0) node[right] {$x$};
\draw[thick,->] (0,-0.2) -- (0,6.1) node[above] {$y$};
\draw (1,-0.2) node[below] {\small 1};
\draw (-0.2,1) node[left] {\small 1};
\draw [gray, thin] (0,1) -- (4.9,5.9);
\draw [very thick,color=blue, name path=sigmarho](-0.2,0.-0.06667)--(5.9,1.96667); 
\draw [very thick,color=blue](-0.2,0.26667)--(5.9,2.3);
\draw [very thick,color=blue](-0.2,0.6)--(5.9,2.63333);

\draw  (0,0) node {$\bullet$};
\draw  (1,1) node {$\bullet$};
\draw  (2,1) node {$\bullet$};
\draw  (2,2) node {$\bullet$};
\draw  (3,1) node {$\bullet$};
\draw  (3,2) node {$\bullet$};
\draw  (3,3) node {$\bullet$};
\draw  (4,1) node {$\bullet$};
\draw  (4,2) node {$\bullet$};
\draw  (4,3) node {$\bullet$};
\draw  (4,4) node {$\bullet$};
\draw  (5,1) node {$\bullet$};
\draw  (5,2) node {$\bullet$};
\draw  (5,3) node {$\bullet$};
\draw  (5,4) node {$\bullet$};
\draw  (5,5) node {$\bullet$};


  \end{scope}
\end{tikzpicture}
\end{minipage}
\hfill
\begin{minipage}{0.39\textwidth}
\begin{tikzpicture}[scale=4]
  \begin{scope}
\draw[ultra thick] (0,0) -- (0,1) -- (1,1) -- (1,0) -- cycle;
\draw[very thick] (0,1) -- (1,0);
\draw[very thick] (0,1) -- (0.5,0);
\draw[very thick] (0.5,1) -- (1,0);

\draw[very thick] (0,1) -- (0.333,0);
\draw[very thick] (0.333,1) -- (0.667,0);
\draw[very thick] (0.667,1) -- (1,0);

\draw[very thick] (0,1) -- (0.25,0);
\draw[very thick] (0.25,1) -- (0.5,0);
\draw[very thick] (0.5,1) -- (0.75,0);
\draw[very thick] (0.75,1) -- (1,0);

\draw[very thick] (0,1) -- (0.2,0);
\draw[very thick] (0.2,1) -- (0.4,0);
\draw[very thick] (0.4,1) -- (0.6,0);
\draw[very thick] (0.6,1) -- (0.8,0);
\draw[very thick] (0.8,1) -- (1,0);



\draw [color=blue] (0.333,0) node {$\bullet$};



\draw [color=blue] (0.33333,0.66667) node {$\bullet$};


\draw [color=blue] (0.33333,0.33333) node {$\bullet$};


  \end{scope}

\end{tikzpicture}
\end{minipage}
\caption{Dualité. Les trois droites de pente 1/3 croisent au total 5+1 points de $P_5$ et correspondent à 3 points de  ${\mathcal A}_5$ qui rencontrent au total 6 droites de $L_5$}\label{f:13}
\end{figure}

Pour trouver la somme \eqref{e1}, classifions les droites $t^*$ selon leur pente. La pente de chaque telle droite est un nombre rationnel $\frac{p}{q}=\frac{l-l'}{k-k'}$, où $0 < p < q \leq n$. Il y a $q$ droites différentes de pente $p/q$ qui passent par des points entiers de l'ensemble $P_n$; elle correspondent à $\rho=0,1/q,\ldots, (q-1)/q$ (voir Fig.~\ref{f:13}). Toutes les valeurs de $\rho$ sont situées dans l'intervalle $[0,1[$, et pour tout $m\in \{0,\ldots,n\}$, il y a exactement une point entier avec $x=m$ croisé par une de ces $q$ droites. Le nombre total de points de l'ensemble $P_n$ qu'elles croisent est donc $n+1$, et leur contribution totale à la somme \eqref{e1} est $n+1-q$.

Le numérateur $p$ de la fraction $p/q$ doit être premier avec le dénominateur $q$. Donc, pour un $q$ donné, il prend exactement $\varphi(q)$ valeurs, où $\varphi$ est la fonction indicatrice d'Euler. Quant à $q$, il prend les valeurs de 2 à $n$. Donc l'équation \eqref{e1} peut être réécrite comme 
\[p_s(n)=1+n+\sum_{q=2}^{n} \varphi(q) (n+1-q) = 1+\sum_{q=1}^{n} \varphi(q) (n+1-q).\]
La dernière égalité est vraie puisque $\varphi(1)=1$. Maintenant, l'estimation $p_s(n)=\frac{n^3}{\pi^2}+O(n^2\log n)$ 
découle immédiatement de la formule connue $\sum_{q=1}^{n} \varphi(q)=\frac{3}{\pi^2}n^2+O(n \log n)$:
\begin{align*}
 \sum_{q=1}^{n} \varphi(q) (n+1-q)&=n\varphi(1)+(n-1)\varphi(2)+\cdots+ \varphi(n)\\
&= \sum_{q=1}^{n} \varphi(q)+\sum_{q=1}^{n-1} \varphi(q)+\cdots \sum_{q=1}^{1} \varphi(q)\\
&= \frac{3}{\pi^2}(n^2+(n-1)^2+\cdots 1^2)+O(n^2\log n)\\
&= \frac{3}{\pi^2} \frac{n^3}{3} +O(n^2\log n),
\end{align*}

 complétant la démonstration du théorème \ref{t:totalSt}. 

\section{Mots de rotation}
\begin{figure}
\centering
\def\ra{0.281966}
\begin{tikzpicture}
  \begin{scope}
\draw[step=1cm,gray,very thin,name path=ggg] (-0.2,-0.2) grid (6.9,3.9);
 \draw[thick,->] (-0.2,0) -- (7.1,0) node[right] {$x$};
\draw[thick,->] (0,-0.2) -- (0,4.1) node[above] {$y$};
\draw (1,-0.2) node[below] {\small 1};
\draw (-0.2,1) node[left] {\small 1};
\draw (-0.2,1-\ra) node[left] {\small $1-\sigma$};

\draw [color=gray,very thin] (-0.2,1-\ra) -- (6.9,1-\ra);
\draw [color=gray,very thin] (-0.2,2-\ra) -- (6.9,2-\ra);
\draw [color=gray,very thin] (-0.2,3-\ra) -- (6.9,3-\ra);
\draw [thick,name path=sigmarho](-0.2,0.381)--(6.9,3.192) node[right] {$y=\alpha x +\rho$};
\path [name path=Y0] (0,0)--(0,5);
\path [name intersections={of=sigmarho and Y0,by=E0}];
\node [fill=black,inner sep=1pt] at (E0) {};

\path [name path=Y1] (1,0)--(1,5);
\path [name intersections={of=sigmarho and Y1,by=E1}];
\node [fill=black,inner sep=1pt] at (E1) {};

\path [name path=Y2] (2,0)--(2,5);
\path [name intersections={of=sigmarho and Y2,by=E2}];
\node [fill=black,inner sep=1pt] at (E2) {};

\path [name path=Y3] (3,0)--(3,5);
\path [name intersections={of=sigmarho and Y3,by=E3}];
\node [fill=black,inner sep=1pt] at (E3) {};

\path [name path=Y4] (4,0)--(4,5);
\path [name intersections={of=sigmarho and Y4,by=E4}];
\node [fill=black,inner sep=1pt] at (E4) {};

\path [name path=Y5] (5,0)--(5,5);
\path [name intersections={of=sigmarho and Y5,by=E5}];
\node [fill=black,inner sep=1pt] at (E5) {};

\path [name path=Y6] (6,0)--(6,5);
\path [name intersections={of=sigmarho and Y6,by=E6}];
\node [fill=black,inner sep=1pt] at (E6) {};

\fill [opacity=0.1](-0.2,1-\ra) -- (-0.2,1) -- (6.9,1) -- (6.9,1-\ra) -- cycle;
\fill [opacity=0.1](-0.2,2-\ra) -- (-0.2,2) -- (6.9,2) -- (6.9,2-\ra) -- cycle;
\fill [opacity=0.1](-0.2,3-\ra) -- (-0.2,3) -- (6.9,3) -- (6.9,3-\ra) -- cycle;

\draw [color=blue] (0, -0.7) node {\it 0};
\draw [color=blue] (1, -0.7) node {\it 1};
\draw [color=blue] (2, -0.7) node {\it 0};
\draw [color=blue] (3, -0.7) node {\it 0};
\draw [color=blue] (4, -0.7) node {\it 0};
\draw [color=blue] (5, -0.7) node {\it 0};
\draw [color=blue] (6, -0.7) node {\it 1};

  \end{scope}
\end{tikzpicture}
\caption{Un mot de rotation $0100001\cdots$}\label{f:rot}
 
\end{figure}
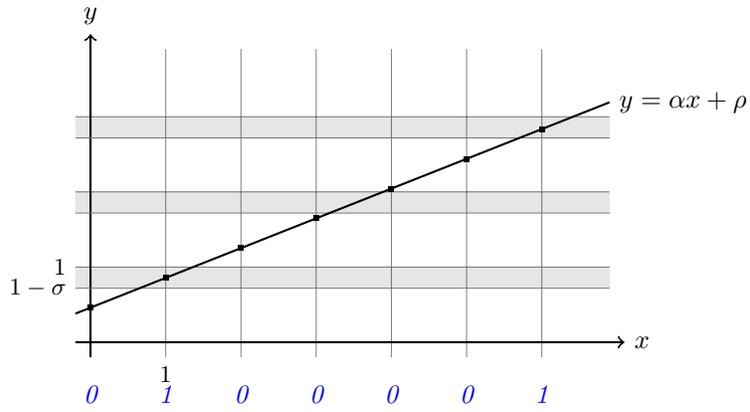

La même méthode de comptage des faces de l'image duale a été utilisée au moins une fois de plus, dans notre article \cite{cf} sur le nombre de facteurs arithmétiques d’un mot sturmien ${\bf s}$, c’est-à-dire la fonction
\[a_{\bf s}(n) = \# \{{\bf s}[k] {\bf s}[k + d] \cdots {\bf s} [k + (n-1) d] | k, d> 0 \}. \]
Pour un mot sturmien ${\bf s} = {\bf s}_{\sigma, \rho}$, le symbole ${\bf s}[k + qd]$ est égal à 1 si et seulement si
\[\lfloor (k+qd+1)\sigma+\rho \rfloor \neq \lfloor (k+qd)\sigma+\rho \rfloor, \]
c'est-à-dire, si
\[\lfloor (qd+1)\sigma+\rho' \rfloor \neq \lfloor qd\sigma+\rho' \rfloor, \]
où $\rho'= \rho + k \sigma$. Cette condition est équivalente à ce que
\[\{qd \sigma + \rho'\}> 1- \sigma. \]
Comme la fonction $a_{{\bf s}_{\sigma, \rho}}(n)$ compte tous les $k$ et $d$ ensemble, et comme la suite $\{d \sigma \}_{d = 1}^\infty$ est uniformément distribuée sur $[0,1]$, alors le nombre de facteurs arithmétiques d’un mot sturmien est égal au nombre total  ${\bf r}_{\alpha,\rho,\sigma}$ de {\it mots de rotation} de longueur $n$ avec le paramètre $\sigma$ fixé. Ici, un mot de rotation  ${\bf r}={\bf r}_{\alpha,\rho,\sigma}$ est défini par
\[{\bf r}[q]=\begin{cases}
              0, \mbox{~si~} \{q\alpha+\rho \}\leq 1-\sigma,\\
              1, \mbox{~sinon}.
             \end{cases}\]
Le paramètre $\sigma$ est donc la largeur des bandes grises à la Fig.~\ref{f:rot} est et fixé, tandis que la pente $\alpha$ et l'intercept $\rho$ de la droite prendent les valeurs dans l'intervalle $[0,1[$.

En particulier, si $\alpha = \sigma$, on a ${\bf r}_{\sigma,\rho,\sigma}={\bf s}_{\sigma,\rho}$, ce qui veut dire qu'un mot sturmien est un cas particulier d'un mot de rotation. 

Dans notre papier avec Julien Cassaigne \cite{cf} notre objectif était donc d'estimer ou de trouver le nombre total $a_{{\bf s}_{\sigma,\rho}}(n)$ de mots de rotation avec le paramètre $\sigma$ fixé; ce nombre sera noté $r_{\sigma}(n)$.

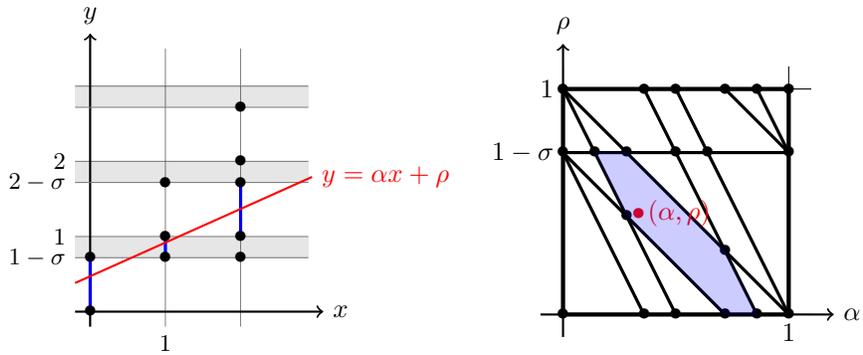
\begin{figure}
 \centering
\begin{minipage}{0.52\textwidth}
 \def\ra{0.281966}
\begin{tikzpicture}
  \begin{scope}
\draw[step=1cm,gray,very thin,name path=ggg] (-0.2,-0.2) grid (2.9,3.5);
 \draw[thick,->] (-0.2,0) -- (3.1,0) node[right] {$x$};
\draw[thick,->] (0,-0.2) -- (0,3.7) node[above] {$y$};
\draw (1,-0.2) node[below] {\small 1};
\draw (-0.2,1) node[left] {\small 1};
\draw (-0.2,1-\ra) node[left] {\small $1-\sigma$};
\draw (-0.2,2) node[left] {\small 2};
\draw (-0.2,2-\ra) node[left] {\small $2-\sigma$};

\draw [color=gray,very thin] (-0.2,1-\ra) -- (2.9,1-\ra);
\draw [color=gray,very thin] (-0.2,2-\ra) -- (2.9,2-\ra);
\draw [color=gray,very thin] (-0.2,3-\ra) -- (2.9,3-\ra);

\fill [opacity=0.1](-0.2,1-\ra) -- (-0.2,1) -- (2.9,1) -- (2.9,1-\ra) -- cycle;
\fill [opacity=0.1](-0.2,2-\ra) -- (-0.2,2) -- (2.9,2) -- (2.9,2-\ra) -- cycle;
\fill [opacity=0.1](-0.2,3-\ra) -- (-0.2,3) -- (2.9,3) -- (2.9,3-\ra) -- cycle;

\draw [very thick, color=blue] (0,0)--(0,1-\ra);
\draw [very thick, color=blue] (1,1-\ra)--(1,1);
\draw [very thick, color=blue] (2,1)--(2,2-\ra);

\draw  (0,0) node {$\bullet$};
\draw  (1,1) node {$\bullet$};
\draw  (2,1) node {$\bullet$};
\draw  (2,2) node {$\bullet$};
\draw (0,1-\ra) node {$\bullet$};
\draw (1,1-\ra) node {$\bullet$};
\draw (2,1-\ra) node {$\bullet$};
\draw (1,2-\ra) node {$\bullet$};
\draw (2,2-\ra) node {$\bullet$};
\draw (2,3-\ra) node {$\bullet$};

\draw [thick,color=red,name path=alpharho](-0.2,0.381)--(2.95,1.792) node[right] {$y=\alpha x +\rho$};

  \end{scope}
\end{tikzpicture}
\end{minipage}
\hfill
\begin{minipage}{0.47\textwidth}
\def\ra{0.281966}
\begin{tikzpicture}[scale=3]
  \begin{scope}
\draw[step=1cm, name path=ggg] (0,0) grid (1.1,1.1);
\draw[thick,->] (-0.1,0) -- (1.2,0) node[right] {$\alpha$};
\draw[thick,->] (0,-0.1) -- (0,1.2) node[above] {$\rho$};
\draw[ultra thick] (0,0) -- (0,1) -- (1,1) -- (1,0) -- cycle;
\draw[very thick] (0,1) -- (1,0);
\draw[very thick] (0,1) -- (0.5,0);
\draw[very thick] (0.5,1) -- (1,0);
\draw[very thick] (0,1-\ra) -- (1, 1-\ra);
\draw[very thick] (0,1-\ra) -- (1-\ra, 0);
\draw[very thick] (1-\ra, 1) -- (1,1-\ra);

\draw[very thick] (0,1-\ra) -- ( 0.5-\ra/2,0);
\draw[very thick] (0.5-\ra/2,1) -- (1-\ra/2, 0);
\draw[very thick] (1-\ra/2, 1) -- (1,1-\ra);

\draw (1,0) node[below] {$1$};
\draw (0,1) node[left] {$1$};
\draw (0,1-\ra) node[left] {$1-\sigma$};

\draw (1, 1) node {$\bullet$};
\draw (1, 0) node {$\bullet$};
\draw (0, 1) node {$\bullet$};
\draw (0, 0) node {$\bullet$};
\draw (0.5, 0) node {$\bullet$};
\draw (0.5, 1) node {$\bullet$};

\draw (1, 1-\ra) node {$\bullet$};
\draw (1-\ra, 0) node {$\bullet$};
\draw (0, 1-\ra) node {$\bullet$};
\draw (1-\ra, 1) node {$\bullet$};
\draw (1-\ra/2, 0) node {$\bullet$};
\draw (1-\ra/2, 1) node {$\bullet$};
\draw (0.5-\ra/2, 0) node {$\bullet$};
\draw (0.5-\ra/2, 1) node {$\bullet$};
\draw (\ra/2,1-\ra) node {$\bullet$};
\draw (\ra,1-\ra) node {$\bullet$};
\draw (0.5+\ra/2,1-\ra) node {$\bullet$};
\draw (0.5,1-\ra) node {$\bullet$};
\draw (1-\ra,\ra) node {$\bullet$};
\draw (\ra,1-2*\ra) node {$\bullet$};

\draw [color=red] (0.3347,0.448) node {$\bullet$} node[right] {$(\alpha,\rho)$};
\fill [color=blue, opacity=0.2] (\ra/2,1-\ra) -- (\ra,1-\ra) -- (1-\ra,\ra) -- (1-\ra/2,0) -- (1-\ra,0) -- (\ra,1-2*\ra) -- cycle;
  \end{scope}

\end{tikzpicture}

\end{minipage}\caption{Dualité pour les mots de rotation avec $\sigma$ fixé}\label{f:d2}
\end{figure}

Comme dans le problème précédent, un mot de rotation est déterminé par la position de la droite $y=\alpha x + \rho$ par rapport aux points entiers $(k,l)$ avec $1\leq l\leq k \leq n$, mais aussi par rapport aux points $(k,l-\sigma)$, où $0\leq l-1\leq k \leq n$ (voir la partie gauche de Fig.~\ref{f:d2}). Sur le graphique dual, ces points correspondent aux droites $k\alpha+\rho=l$ et $k\alpha+\rho=l-\sigma$ (voir la partie droite de la même figure). La  différence principale par rapport au problème précédent est que l'arrangement d’ordre $n$ correspond aux mots de longueur $n + 1$, car dans ce problème, la longueur d'un mot est le nombre de lignes entières verticales croisées, $x=0$ compris, et pas le nombre d'intervalles entre deux telles lignes. Fig.~\ref{f:d2} montre un arrangement d’ordre 2 correspondant aux mots de longueur 3; la face coloriée correspond à la ligne passant par les intervalles verticaux $y\in ]0,1-\sigma[$ pour $x=0$, $y\in ]1-\sigma,1[$ pour $x=1$ et $y\in ]1,2-\sigma[$ pour $x=2$. Le 
mot correspondant est donc $010$.

Comme dans le problème précédent, chaque face de l'arrangement de la Fig.~\ref{f:d2} correspond à un mot de rotation, et donc  le nombre de faces $f_{\sigma}(n)$ de l'arrangement d'ordre $n$ est une borne supérieure pour $r_{\sigma}(n)$. Grâce aux les mêmes arguments que dans le papier de Berstel et Pocchiola \cite{bp}, nous obtenons que, à condition que $\sigma$ soit irrationnel, le nombre de faces $f(n) = f_{\sigma}(n)$ ne dépend pas de $\sigma$ et est égal à

\[f(n)=2+\frac{n(n+1)(n+2)}{3}+2 \sum_{q=1}^n (n-q+1) \varphi(q).\]

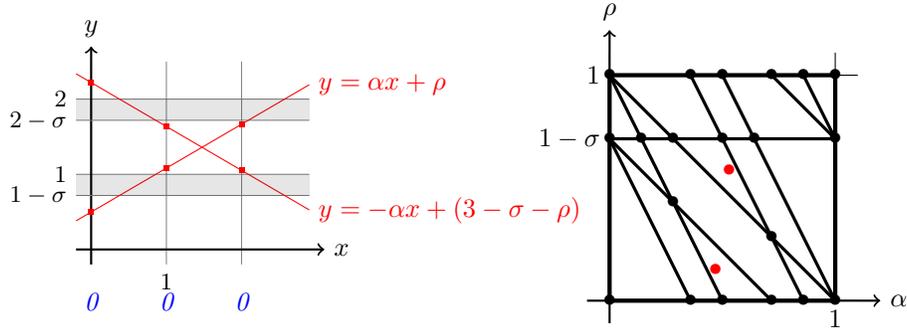
\begin{figure}
\centering
\begin{minipage}{0.54\textwidth}
\def\ra{0.281966}
\begin{tikzpicture}
  \begin{scope}
\draw[step=1cm,gray,very thin,name path=ggg] (-0.2,-0.2) grid (2.9,2.5);
 \draw[thick,->] (-0.2,0) -- (3.1,0) node[right] {$x$};
\draw[thick,->] (0,-0.2) -- (0,2.7) node[above] {$y$};
\draw (1,-0.2) node[below] {\small 1};
\draw (-0.2,1) node[left] {\small 1};
\draw (-0.2,1-\ra) node[left] {\small $1-\sigma$};
\draw (-0.2,2) node[left] {\small 2};
\draw (-0.2,2-\ra) node[left] {\small $2-\sigma$};

\draw [color=gray,very thin] (-0.2,1-\ra) -- (2.9,1-\ra);
\draw [color=gray,very thin] (-0.2,2-\ra) -- (2.9,2-\ra);
\draw [color=red,name path=sigmarho](-0.2,0.381)--(2.9,2.192) node[right] {$y=\alpha x +\rho$};

\draw [color=red,name path=sigmarho2](-0.2,3-\ra-0.381)--(2.9,3-\ra-2.192) node[right] {$y=-\alpha x +(3-\sigma-\rho)$};
\path [name path=Y0] (0,0)--(0,3);
\path [name intersections={of=sigmarho and Y0,by=E0}];
\path [name intersections={of=sigmarho2 and Y0,by=F0}];
\node [fill=red,inner sep=1pt] at (E0) {};
\node [fill=red,inner sep=1pt] at (F0) {};

\path [name path=Y1] (1,0)--(1,3);
\path [name intersections={of=sigmarho and Y1,by=E1}];
\path [name intersections={of=sigmarho2 and Y1,by=F1}];
\node [fill=red,inner sep=1pt] at (F1) {};
\node [fill=red,inner sep=1pt] at (E1) {};

\path [name path=Y2] (2,0)--(2,3);
\path [name intersections={of=sigmarho and Y2,by=E2}];
\path [name intersections={of=sigmarho2 and Y2,by=F2}];
\node [fill=red,inner sep=1pt] at (F2) {};
\node [fill=red,inner sep=1pt] at (E2) {};

\fill [opacity=0.1](-0.2,1-\ra) -- (-0.2,1) -- (2.9,1) -- (2.9,1-\ra) -- cycle;
\fill [opacity=0.1](-0.2,2-\ra) -- (-0.2,2) -- (2.9,2) -- (2.9,2-\ra) -- cycle;

\draw [color=blue] (0, -0.7) node {\it 0};
\draw [color=blue] (1, -0.7) node {\it 0};
\draw [color=blue] (2, -0.7) node {\it 0};

  \end{scope}
\end{tikzpicture}
\end{minipage}
\hfill
\begin{minipage}{0.42\textwidth}
\def\ra{0.281966}
\begin{tikzpicture}[scale=3]
\draw[step=1cm, name path=ggg] (0,0) grid (1.1,1.1);
\draw[thick,->] (-0.1,0) -- (1.2,0) node[right] {$\alpha$};
\draw[thick,->] (0,-0.1) -- (0,1.2) node[above] {$\rho$};
\draw[ultra thick] (0,0) -- (0,1) -- (1,1) -- (1,0) -- cycle;
\draw[very thick] (0,1) -- (1,0);
\draw[very thick] (0,1) -- (0.5,0);
\draw[very thick] (0.5,1) -- (1,0);
\draw[very thick] (0,1-\ra) -- (1, 1-\ra);
\draw[very thick] (0,1-\ra) -- (1-\ra, 0);
\draw[very thick] (1-\ra, 1) -- (1,1-\ra);

\draw[very thick] (0,1-\ra) -- ( 0.5-\ra/2,0);
\draw[very thick] (0.5-\ra/2,1) -- (1-\ra/2, 0);
\draw[very thick] (1-\ra/2, 1) -- (1,1-\ra);

\draw (1,0) node[below] {$1$};
\draw (0,1) node[left] {$1$};
\draw (0,1-\ra) node[left] {$1-\sigma$};

\draw (1, 1) node {$\bullet$};
\draw (1, 0) node {$\bullet$};
\draw (0, 1) node {$\bullet$};
\draw (0, 0) node {$\bullet$};
\draw (0.5, 0) node {$\bullet$};
\draw (0.5, 1) node {$\bullet$};

\draw (1, 1-\ra) node {$\bullet$};
\draw (1-\ra, 0) node {$\bullet$};
\draw (0, 1-\ra) node {$\bullet$};
\draw (1-\ra, 1) node {$\bullet$};
\draw (1-\ra/2, 0) node {$\bullet$};
\draw (1-\ra/2, 1) node {$\bullet$};
\draw (0.5-\ra/2, 0) node {$\bullet$};
\draw (0.5-\ra/2, 1) node {$\bullet$};
\draw (\ra/2,1-\ra) node {$\bullet$};
\draw (\ra,1-\ra) node {$\bullet$};
\draw (0.5+\ra/2,1-\ra) node {$\bullet$};
\draw (0.5,1-\ra) node {$\bullet$};
\draw (1-\ra,\ra) node {$\bullet$};
\draw (\ra,1-2*\ra) node {$\bullet$};

\draw [color=red] (0.53, 0.58) node {$\bullet$};
\draw [color=red] (1-0.53, 1-\ra-0.58) node {$\bullet$};

\end{tikzpicture}
\end{minipage}
\caption{Droites et points symétriques}\label{f:sym}
\end{figure}

Cependant, contrairement au cas précédent, il n'est pas vrai que des faces différentes correspondent à des mots de rotation différents! En particulier, pour chaque $\alpha$ et $\rho$ on a
\begin{equation}\label{e:sym}
{\bf r}_{\alpha,\rho,\sigma}={\bf r}_{\{-\alpha\},\{-\sigma-\rho\},\sigma}
\end{equation}
en raison de la symétrie représentée à gauche de la Fig.~\ref{f:sym}: on voit que les droites $y =\alpha x + \rho$ et 
$y = -\alpha x + (k-\sigma-\rho)$ pour $ k \in \mathbb Z$ se comportent symétriquement par rapport aux bandes 
$\{y \} <1- \sigma$; comme on peut ajouter un entier à l'un des coefficients sans modifier la partie fractionnaire, cela justifie l'égalité \eqref{e:sym}. Ainsi, le graphique dual est également symétrique avec deux point symétriques pointés en rouge (voir la partie droite de la Fig.~\ref{f:sym}), ce qui signifie immédiatement que
\[r_{\sigma}(n+1)\leq \frac{f(n)}{2}+1,\]
car les centres de symétrie sont situés à l'intérieur de deux faces.

En outre, il existe également d'autres faces correspondant au même mot! Pour les classer, nous avons dû faire une expérience de calcul avec une variante de la technique de Monte-Carlo. Nous avons trouvé un point $(\alpha,\rho)$ sur chaque face de l'arrangement d'un ordre donné et nous avons étudié ce qui se passe exactement sur les faces non symétriques qui continuent à donner le même mot de rotation. Le résultat dépend de $\sigma$, et nous n’avons obtenu des formules précises que pour certaines de ses valeurs, à savoir pour $1/3 <\sigma <2/3$. Par exemple, si 
 $3/8<\sigma<2/5$, on a
\[r_{\sigma}(n+1)=\begin{cases}
                 \frac{f(n)}{2}-7, \mbox{~si~} n \mbox{~est pair},\\
		  \frac{f(n)}{2}-8, \mbox{~si~} n \mbox{~est impair}
                \end{cases}\]
à partir de $n=8$.

Pour $1/3<\sigma<2/3$, la différence entre $r_{\sigma}(n+1)$ et $\frac{f(n)}{2}$ est bornée, et on sait donc que l'ordre de croissance de $r_{\sigma}(n)$ est $n^3(1/6+1/\pi^2)$. Pour les autres valeurs de $\sigma$, nous n'avons que des bornes supérieures et inférieures pour $r_{\sigma}(n)$, toutes d'ordre $n^3$ \cite{cf,f_lower}.

La prochaine question raisonnable est la suivante: quel est le nombre total de mots de rotation, toutes valeurs de $\sigma$, $\rho$ et $\alpha$ prises ensemble? L’idée de considérer un arrangement tridimensionnel et de compter les zones coupées par des plans semble attrayante, mais c’est en fait une autre technique qui a été utilisée pour trouver cette valeur dans notre papier avec D. Jamet \cite{fj}.

À la fin de cette partie du texte, je souligne que les expériences de calcul sont souvent utilisées dans ce domaine pour compter les mots jusqu'à une certaine longueur. Cela permet de formuler la conjecture et d'initialiser la démonstration par récurrence. Les mathématiques sont devenues des sciences expérimentales et, dans ce sens, plus proches que jamais de l'informatique théorique.

\section{Palindromes dans les mots sturmiens}
Comme d'habitude, un {\it palindrome} est un mot fini $u$ qui se lit de la même manière de gauche à droite et de droite à gauche, c'est-à-dire, $u = u[1] \cdots u[n] = u[n] \cdots u[1]$, comme dans {\it kayak} ou {\it 01001010010}. Nous disons aussi qu'un palindrome $aua$, où $a$ est une lettre, est une {\it extension palindromique} du palindrome $u$.

Tout mot sturmien est {\it riche} de palindromes dans un sens mathématique strict que nous allons définir. Mais la première proposition concerne le nombre de facteurs palindromiques d'une longueur donnée dans un mot sturmien.

\begin{proposition}[\cite{droubay_pirillo}]
 Pour tout mot sturmien ${\bf s}$, le nombre $h_{\bf s}(n)$ de ses facteurs palindromes de longueur $n$ est $2$ si $n$ est impair et $1$ si $n$ est pair.
\end{proposition}
\noindent {\sc Démonstration.} {\'E}videmment, ${\bf s}$ contient les lettres $0$ et $1$ et exactement l'un des mots $00$ ou $11$. Ainsi, l'énoncé est vrai pour $n = 1$ et $n = 2$. Notons maintenant qu'un facteur palindrome $u'= {\bf s}[i..n + i + 1]$ de longueur $n$ de ${\bf s}$ est une extension palindromique du palindrome $u = {\bf s}[i + 1..n + i]$ de longueur $n$ qui est également un facteur de ${\bf s}$. De plus, pour chaque palindrome $u$, au plus une de ses extensions $0u0$ et $1u1$ peut être facteur du même mot sturmien, car un mot sturmien est équilibré et $|0u0|_0-|1u1|_0=2>1$. Nous avons donc $h_{\bf s}(n+2)\leq h_{\bf s}(n)$. Il reste à trouver dans ${\bf s}$ deux palindromes de toute longueur impaire et un  palindrome de toute longueur paire.

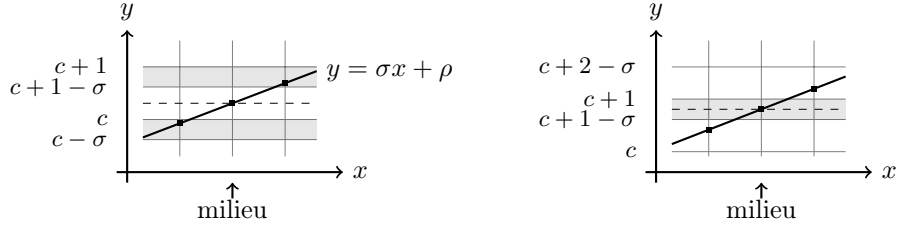
\begin{figure}
\begin{minipage}{0.57\textwidth}
\def\ra{0.381966}
\begin{tikzpicture}[scale=.7]
  \begin{scope}
\draw[step=1cm,gray,very thin,name path=ggg] (0.3,0.3) grid (3.6,2.5);
 \draw[thick,->] (-0.2,0) -- (4.1,0) node[right] {$x$};
\draw[thick,->] (0,-0.2) -- (0,2.7) node[above] {$y$};
\draw (-0.2,1) node[left] {\small $c$};
\draw (-0.2,1-\ra) node[left] {\small $c-\sigma$};
\draw (-0.2,2) node[left] {\small $c+1$};
\draw (-0.2,2-\ra) node[left] {\small $c+1-\sigma$};

\draw [color=gray,very thin] (0.3,1-\ra) -- (3.6,1-\ra);
\draw [color=gray,very thin] (0.3,2-\ra) -- (3.6,2-\ra);

\draw [dashed,very thin] (0.3,1.5-\ra/2) -- (3.6,1.5-\ra/2);

\def\rrho{1.5-2.5*\ra}
\draw [thick,name path=sigmarho](0.3,\ra*0.3+\rrho)--(3.6,\ra*3.6+\rrho) node[right] {$y=\sigma x +\rho$};

\node [fill,inner sep=1pt] at (2,1.5-\ra/2) {};
\node [fill,inner sep=1pt] at (1,1.5-\ra/2-\ra) {};
\node [fill,inner sep=1pt] at (3,1.5-\ra/2+\ra) {};

\fill [opacity=0.1](0.3,1-\ra) -- (0.3,1) -- (3.6,1) -- (3.6,1-\ra) -- cycle;
\fill [opacity=0.1](0.3,2-\ra) -- (0.3,2) -- (3.6,2) -- (3.6,2-\ra) -- cycle;

\draw[thick,->] (2,-0.5) -- (2,-0.2);
\draw (2, -0.7) node {milieu};

  \end{scope}
 \end{tikzpicture}
\end{minipage}
\hfill
\begin{minipage}{0.42\textwidth}
\def\ra{0.611}
\begin{tikzpicture}[scale=.7]
  \begin{scope}
\draw[step=1cm,gray,very thin,name path=ggg] (0.3,0.3) grid (3.6,2.5);
 \draw[thick,->] (-0.2,0) -- (4.1,0) node[right] {$x$};
\draw[thick,->] (0,-0.2) -- (0,2.7) node[above] {$y$};
\draw (-0.2,1) node[left] {\small $c+1-\sigma$};
\draw (-0.2,1-\ra) node[left] {\small $c$};
\draw (-0.2,2) node[left] {\small $c+2-\sigma$};
\draw (-0.2,2-\ra) node[left] {\small $c+1$};

\draw [color=gray,very thin] (0.3,1-\ra) -- (3.6,1-\ra);
\draw [color=gray,very thin] (0.3,2-\ra) -- (3.6,2-\ra);

\draw [dashed,very thin] (0.3,1.5-\ra/2) -- (3.6,1.5-\ra/2);

\def\rrho{1.5-2.5-2.5*\ra}
\draw [thick,name path=sigmarho](0.3,1.5-\ra/2-1.7+1.7*\ra)--(3.6,1.5-\ra/2+1.6-1.6*\ra); 

\node [fill,inner sep=1pt] at (2,1.5-\ra/2) {};
\node [fill,inner sep=1pt] at (1,1.5-\ra/2-1+\ra) {};
\node [fill,inner sep=1pt] at (3,1.5-\ra/2+1-\ra) {};

\fill [opacity=0.1](0.3,1) -- (0.3,2-\ra) -- (3.6,2-\ra) -- (3.6,1) -- cycle;

\draw[thick,->] (2,-0.5) -- (2,-0.2);
\draw (2, -0.7) node {milieu};

  \end{scope}
\end{tikzpicture}
 \end{minipage}
\caption{Palindromes sturmiens de longueur impaire}\label{f:pal1}

\end{figure}

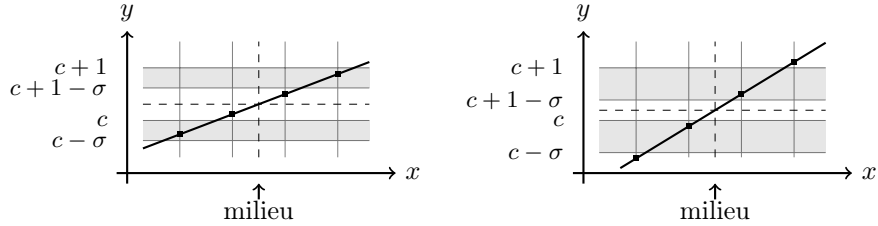
\begin{figure}
\centering
 \begin{minipage}{0.49\textwidth}
\def\ra{0.381966}
\begin{tikzpicture}[scale=.7]
  \begin{scope}
\draw[step=1cm,gray,very thin,name path=ggg] (0.3,0.3) grid (4.6,2.5);
 \draw[thick,->] (-0.2,0) -- (5.1,0) node[right] {$x$};
\draw[thick,->] (0,-0.2) -- (0,2.7) node[above] {$y$};
\draw (-0.2,1) node[left] {\small $c$};
\draw (-0.2,1-\ra) node[left] {\small $c-\sigma$};
\draw (-0.2,2) node[left] {\small $c+1$};
\draw (-0.2,2-\ra) node[left] {\small $c+1-\sigma$};

\draw [color=gray,very thin] (0.3,1-\ra) -- (4.6,1-\ra);
\draw [color=gray,very thin] (0.3,2-\ra) -- (4.6,2-\ra);

\draw [dashed,very thin] (0.3,1.5-\ra/2) -- (4.6,1.5-\ra/2);
\draw [dashed,very thin] (2.5,0.3) -- (2.5,2.5);

\def\rrho{1.5-3*\ra}
\draw [thick,name path=sigmarho](0.3,\ra*0.3+\rrho)--(4.6,\ra*4.6+\rrho);

\node [fill,inner sep=1pt] at (2,1.5-\ra) {};
\node [fill,inner sep=1pt] at (1,1.5-2*\ra) {};
\node [fill,inner sep=1pt] at (3,1.5) {};
\node [fill,inner sep=1pt] at (4,1.5+\ra) {};

\fill [opacity=0.1](0.3,1-\ra) -- (0.3,1) -- (4.6,1) -- (4.6,1-\ra) -- cycle;
\fill [opacity=0.1](0.3,2-\ra) -- (0.3,2) -- (4.6,2) -- (4.6,2-\ra) -- cycle;

\draw[thick,->] (2.5,-0.5) -- (2.5,-0.2);
\draw (2.5, -0.7) node {milieu};

  \end{scope}
 \end{tikzpicture}  
 \end{minipage}
\hfill
 \begin{minipage}{0.49\textwidth}
\def\ra{0.611}
\begin{tikzpicture}[scale=.7]
  \begin{scope}
\draw[step=1cm,gray,very thin,name path=ggg] (0.3,0.3) grid (4.6,2.5);
 \draw[thick,->] (-0.2,0) -- (5.1,0) node[right] {$x$};
\draw[thick,->] (0,-0.2) -- (0,2.7) node[above] {$y$};
\draw (-0.2,1) node[left] {\small $c$};
\draw (-0.2,1-\ra) node[left] {\small $c-\sigma$};
\draw (-0.2,2) node[left] {\small $c+1$};
\draw (-0.2,2-\ra) node[left] {\small $c+1-\sigma$};

\draw [color=gray,very thin] (0.3,1-\ra) -- (4.6,1-\ra);
\draw [color=gray,very thin] (0.3,2-\ra) -- (4.6,2-\ra);

\draw [dashed,very thin] (0.3,1.5-\ra/2) -- (4.6,1.5-\ra/2);
\draw [dashed,very thin] (2.5,0.3) -- (2.5,2.5);

\def\rrho{1.5-3*\ra}
\draw [thick,name path=sigmarho](0.7,\ra*0.7+\rrho)--(4.6,\ra*4.6+\rrho); 

\node [fill,inner sep=1pt] at (2,1.5-\ra) {};
\node [fill,inner sep=1pt] at (1,1.5-2*\ra) {};
\node [fill,inner sep=1pt] at (3,1.5) {};
\node [fill,inner sep=1pt] at (4,1.5+\ra) {};

\fill [opacity=0.1](0.3,1-\ra) -- (0.3,1) -- (4.6,1) -- (4.6,1-\ra) -- cycle;
\fill [opacity=0.1](0.3,2-\ra) -- (0.3,2) -- (4.6,2) -- (4.6,2-\ra) -- cycle;

\draw[thick,->] (2.5,-0.5) -- (2.5,-0.2);
\draw (2.5, -0.7) node {milieu};

  \end{scope}
 \end{tikzpicture} 
 \end{minipage}
\hfill
\caption{Palindromes sturmiens de longueur paire, $\sigma<1/2$ (à gauche) et $\sigma>1/2$ (à droite)}\label{f:pal2}
\end{figure}

Une des manières les plus simples de les voir est d'utiliser la définition \eqref{e:frpart} et le fait que les facteurs d'un mot sturmien de pente $\sigma$ sont exactement les mots correspondant à toutes les valeurs possibles de $\rho$. Pour trouver un facteur palindrome $u = {\bf s}_{\sigma,\rho}[1..n]$, nous indiquons simplement des valeurs appropriées de $\rho$. Pour $n = 2k + 1$, le 
$\rho$ est choisi de telle sorte que $\{(k+1)\sigma+\rho\}=(1-\sigma)/2$ pour le facteur sturmien avec un 0 au milieu et $\{(k+1)\sigma+\rho\}=1-\sigma/2$ pour le facteur avec un 1 au milieu. Ce choix veut dire que le milieu du mot est situé exactement au centre d'une bande blanche ou grise (voir Fig.~\ref{f:pal1}), et donc les deux morceaux de droites sont symétriques par rapport aux bandes horizontales grises. Les deux facteurs respectifs sont donc des palindromes. Pour une longueur paire $n = 2k$, la situation où $\{k\sigma+\rho\}=1-\sigma$ et donc $\{(k+1)\sigma+\rho\}=0$  est symétrique par rapport aux mêmes bandes mais ne donne pas de palindrome, car les bornes de la bande de largeur $\sigma$ donnent deux symboles différents. Le seul palindrome est donc donné par la valeur de $\rho$ telle que
$\{k\sigma+\rho\}=\{1/2-\sigma\}$ et donc $\{(k+1)\sigma+\rho\}=1/2$: le milieu de la droite est situé au milieu d'une bande blanche, comme il est illustrée à la Fig.~\ref{f:pal2}. Cela prouve l'existence des palindromes mentionnés dans l'énoncé de la proposition.
 \hfill $\Box$

Droubay et Pirillo \cite{droubay_pirillo} ont également montré que cette suite périodique $1,2,1,2,\ldots$ pour le nombre $h_{\bf s}(n)$ de facteurs palindromiques caractérise les mots sturmiens. Ses valeurs ne sont en aucune façon petites. Tout d'abord, il est facile de définir un mot infini sans longs palindromes, comme par exemple $abcabcabc\cdots$. En plus, il est connu que la fonction $h(n)$ est bornée pour tout mot infini de complexité linéaire \cite{abcd}. Enfin, tout mot sturmien est {\it riche} en palindromes de la manière que nous définissons précisément ci-dessous.

Notons le nombre de facteurs palindromiques différents dans un mot fini $u$ par $P(u)$.
\begin{proposition}{\cite{djp}}\label{p:uux}
Pour tout mot fini $u$ et pour toute lettre $x$, on a $P(ux)\leq P(u)+1$.
\end{proposition}
\noindent {\sc Démonstration.} Considérons tous les suffixes de $ux$ qui sont palindromes. S'il n'y en a pas, on a $P(ux)=P(u)$. S'ils existent, ce n'est que le plus long d'entre eux qui peut apparaître dans $ux$ pour la première fois. Effectivement, si $|p'|<|p|$ pour deux tels suffixes $p$ et $p'$ de $ux$, alors $p'$ est aussi un préfixe de $p$, et donc il serait déjà apparu dans $u$. Donc $P(ux)\leq P(u)+1$. \hfill $\Box$

\begin{corollaire}{\cite{djp}}
 Pour tout mot fini $u$, on a $P(u)\leq |u|+1$.
\end{corollaire}
\noindent {\sc Démonstration.} Ceci résulte immédiatement de la proposition \ref{p:uux} et du fait que le mot vide $\epsilon$ est un palindrome, et donc $P(\epsilon)=1$. \hfill $\Box$

On dit qu'un mot fini $u$ avec $P(u)= |u|+1$ est {\it riche} en palindromes.

\begin{proposition}{\cite{djp}}
L'ensemble des mots riches est factoriel, c'est-à-dire, si un mot fini est riche, alors tous ses facteurs sont riches aussi. 
\end{proposition}
\noindent {\sc Démonstration.} Selon la proposition \ref{p:uux}, un mot ne peut être riche que si son préfixe obtenu  en effaçant sa dernière lettre est riche.  Symétriquement, ceci est également vrai pour les suffixes: un mot ne peut être riche que si son suffixe obtenu en effaçant sa première lettre est riche. La proposition suit par récurrence. \hfill $\Box$

On dit qu'un mot infini est riche si tous ses facteurs sont riches.

Les mots sturmiens sont une famille classique de mots riches \cite{djp}. Cependant, comme beaucoup d'autres faits, nous ne le démontrons pas ici.

\section{Suites directrices et système de numération d'Ostrowski}\label{s:directive}

Considérons une {\it suite directrice} ${\bf d}=(d_0,d_1,\ldots,d_n,\ldots)$, où les $d_i$ sont des nombres entiers, $d_i\geq 1$ pour $i>0$ et $d_0\geq 0$. La {\it suite standard} correspondante $(s_n)$ de mots sur l'alphabet $\{a,b\}$ est définie comme suit : 
\begin{equation}\label{e:def}
s_{-1}=b, s_0=a, s_{n+1}=s_n^{d_n}s_{n-1} \mbox{~pour tout~} n\geq 0.
\end{equation}
Le mot $s_n$ s'appelle aussi le  {\it mot standard} d'ordre $n$.

Notons que pour obtenir tous les mots standards possibles, nous devons permettre $d_0 = 0$, mais en raison de la symétrie entre $a$ et $b$, nous pouvons nous limiter au cas de $d_0> 0$.

Il est facile de voir qu’à partir de $n = 0$, chaque $s_n$ est un préfixe de $s_{n + 1}$, et que les longueurs de $s_{n}$ croissent strictement, de sorte qu’il existe un mot infini à droite ${\bf s_d} = \lim_{n \to \infty} s_n$. Le mot 
${\bf s_d}$ est appelé un {\it mot caractéristique} associé à la suite ${\bf d}$. Notons que nous avons renommé l’alphabet de $\{0,1 \}$ en $\{a, b \}$ pour des raisons de commodité, en utilisant le codage $t: 0 \to a, 1 \to b$.

Le fait que le mot caractéristique construit avec \eqref{e:def} est sturmien au sens de \eqref{e:frpart} a été prouvé dans \cite{fmt}; voir aussi \cite{lothaire} pour les démonstrations et les discussions. La pente $\sigma$ de ${\bf s_d}$ est le nombre 
tel que sa fraction continue est
\[\sigma=[0 ; 1 + d_0 ; d_1 ; d_2 ; \ldots ],\]
c'est-à-dire,
\[\sigma=\cfrac{1}{d_0+1+\cfrac{1}{d_1+\cfrac{1}{d_2+\cdots}}}.\]

En plus, ${\bf s_d}=t({\bf s}_{\sigma,\sigma})$, ce qui signifie que ${\bf s_d}$, avec $a$ et $b$ renommés en $0$ et $1$, est obtenu à partir du mot ${\bf s}_{\sigma,0}$ en supprimant le premier symbole (dont l'indice est $0$). Notons aussi que  le mot mécanique inférieur ${\bf s}_{\sigma,0}$ et le mot mécanique supérieur ${\bf s'}_{\sigma,0}$ ne diffèrent qu'au symbole numéro 0 (qui est égal à 0 ou 1, respectivement). Cela veut dire que les deux mots $a{\bf s_d}=t(0s_{\sigma,\sigma})$ et $b{\bf s_d}=t(1s_{\sigma,\sigma})$ sont sturmiens avec le même ensemble de facteurs Fac$(\sigma)$, et donc tout préfixe du mot caractéristique ${\bf s_d}$ est spécial à gauche dans Fac$(\sigma)$.

\begin{exemple}[Mot de Fibonacci]\label{ex:fib1}
 {\rm La suite directrice $(1,1,1,\ldots)$ engendre les mots standards
\[s_{-1}=b, s_0=a, s_1=ab, s_2=aba, s_3=abaab, s_4=abaababa,\]
et ainsi de suite; on voit que par construction, la longueur de chaque $s_n$ est un nombre de Fibonacci, et il est donc raisonnable d'appeler le mot sturmien caractéristique correspondant
\[{\bf f}=abaababaabaababaababa\cdots\]
le {\it mot de Fibonacci }. Sa pente est égale à
\[\tau=\cfrac{1}{2+\cfrac{1}{1+\cfrac{1}{1+\cdots}}}=\frac{3-\sqrt{5}}{2}=0.381966\cdots.\]
}
\end{exemple}

\medskip
Dans la section précédente, nous avons numéroté les symboles dans un mot infini à partir de 0; dans cette section, il est beaucoup plus pratique de commencer les indices à partir de 1 et de considérer un mot sturmien caractéristique
\[{\bf s_d}={\bf s}[1]{\bf s}[2]\cdots {\bf s}[n] \cdots,\]
où
\[{\bf s}[n]=\lfloor (n+1)\sigma \rfloor - \lfloor n\sigma \rfloor.\]
On dit donc que le mot sturmien caractéristique ${\bf s_d}$ est égal à ${\bf s}_{\sigma,0}$, sans oublie de commencer les indices à partir de 1.

Rappelons que, étant donné un mot infini $w$, un palindrome $w(p_1-1..p_2 + 1]$ est appelé l'{\it extension  palindromique} du mot $w(p_1..p_2]$, qui est également un palindrome. Nous pouvons continuer à l'étendre jusqu'à ce que soit nous obtenions un préfixe palindrome $w(0..p_2 + p_1]$ de $w$ soit le mot $w(p_1-d-1..p_2 + d + 1]$ n'est pas un palindrome pour un certain $d <p_1$. Dans les deux cas, on parle de l'{\it extension (palindromique) maximale } $w(0..p_2 + p_1]$ ou $w(p_1-d. .p_2 + d]$ de l'occurrence $w(p_1..p_2]$.

L'ensemble des facteurs d'un mot sturmien est fermé par image miroir, ce qui est facile à voir à partir des Fig.~\ref{f:pal1} et \ref{f:pal2}. Les mots spéciaux à gauche d’un language sturmien Fac$(\sigma)$ sont les préfixes de ${\bf s_d} = {\bf s}_{\sigma, 0}$, et il n’existe aucun autre mot spécial à gauche dans Fac$(\sigma)$, puisque le nombre d'éléments de longueur $n + 1$ dans Fac$(\sigma)$ est le nombre d'éléments de longueur $n$ plus 1. Donc, symétriquement, les mots spéciaux à droite de Fac$(\sigma)$ sont les images miroir des facteurs spéciaux à gauche. Cela signifie immédiatement que les facteurs bispéciaux d'un mot sturmien caractéristique sont exactement ses préfixes qui sont des palindromes.

\begin{proposition}\label{p:bisp}
Pour toute occurrence d'un palindrome dans un mot sturmien caractéristique  ${\bf s}$, son extension palindromique maximale est bispéciale. 
\end{proposition}
 \noindent {\sc Démonstration.} Si l'extension palindromique maximale  est un préfixe de ${\bf s}$, elle est bispéciale comme discuté ci-dessus. Si elle est de la forme ${\bf s}(p_1-d..p_2+d]$ avec $d<p_1$, où ${\bf s}(p_1..p_2]$ est le palindrome initial, cela veut dire que les lettres qui la continuent à droite et à gauche sont différentes : ${\bf s}[p_1-d] \neq {\bf s}[p_2+d+1]$. Comme l'ensemble des facteurs de ${\bf s}$ est fermé par image miroir, et comme le mot ${\bf s}(p_1-d..p_2+d]$ est un palindrome, cela veut dire qu'il est bispécial.  \hfill $\Box$

Les facteurs bispéciaux d'un mot caractéristique construit avec une suite directrice donnée, qui sont aussi exactement les préfixes du mot caractéristique qui sont des palindromes, sont complètement décrits (voir \cite{dlm,lothaire}). Ils sont exactement les mots {\it centraux} $c_{n,j}$, où $0\leq j \leq d_n$, définis comme suit: 
\begin{itemize}
 \item 
Pour tout $n\geq 0$, le mot $c_n$ est obtenu du mot $s_ns_{n-1}$ en effaçant ses deux dernières lettres;
\item
Pour tout $j\geq 0$,  où
 $0\leq j\leq d_n$, le mot $c_{n,j}$ est défini comme $c_{n,j}=s_n^{j}c_{n}$, c'est-à-dire, $c_{n,j}$ est le mot $s_n^{j+1} s_{n-1}$ sans ses deux dernières lettres.
\end{itemize}

On a notamment $c_n=c_{n,0}$ pour tout $n$. Le mot $c_0$ est le mot vide, et on a $c_{0,j}=a^j$, $c_1=a^{d_0}=c_{0,d_{0}}$, $c_{1,j}=(a^{d_0}b)^ja^{d_0}$, $c_2=c_{1,d_1}$, et ainsi de suite; pour tout $n$, on a $c_{n+1}=c_{n,d_n}$.

Notons  \cite{dlm,lothaire} que pour tout
 $n\geq 0$, on a $s_ns_{n-1}\neq s_{n-1}s_n$, mais ces deux mots coïncident à l'exception des deux derniers symboles: $s_n s_{n-1}=c_{n}xy$ et $s_{n-1}s_n=c_{n} yx$, où si $n$ est pair, $x=a$ et $y=b$, et si $n$ est impair, $x=b$ et $y=a$. Le mot central $c_{n}$ s'obtient donc en effaçant les deux derniers symboles de $s_{n-1}s_n$ ou bien de $s_{n}s_{n-1}$.

\begin{exemple}
 {\rm
Les facteurs bispéciaux, ou les mots centraux, ou, ce qui est la même chose, les préfixes palindromes du mot de Fibonacci de l'exemple \ref{ex:fib1} sont
\begin{align*}
 c_0=c_{0,0}&=\epsilon,\\
c_1=c_{1,0}=c_{0,1}&=a,\\
c_2=c_{2,0}=c_{1,1}&=aba,\\
c_3=c_{3,0}=c_{2,1}&=abaaba,\\
c_4=c_{4,0}=c_{3,1}&=abaababaaba,
\end{align*}
etc; tout $c_i$ s'obtient de $s_i$ en effaçant les deux dernières lettres.

Dans le mot sturmien caractéristique $aabaabaaabaabaaaabaaba\cdots$ qui correspond à la suite directrice  ${\bf d_2}=(2,2,2,\ldots)$, on a $s_0=a$, $s_1=aab$, $s_2=aabaaba$, et ainsi de suite. Les facteurs bispéciaux, ou les mots centraux, ou encore les préfixes palindromes de ce mot sont
\begin{align*}
 c_0=c_{0,0}&=\epsilon,\\
c_{0,1}&=a,\\
c_1=c_{1,0}=c_{0,2}&=aa,\\
c_{1,1}&=aabaa,\\
c_2=c_{2,0}=c_{1,2}&=aabaabaa,\\
c_{2,1}&=aabaabaaabaabaa,\\
c_3=c_{3,0}=c_{2,2}&=aabaabaaabaabaaabaabaa,
\end{align*}
etc.
}
\end{exemple}

Le fait suivant est une reformulation de résultats de \cite{dlm,berstel1999}.
\begin{proposition}\label{p:central}
Les facteurs bispéciaux d'un mot sturmien caractéristique ${\bf s_d}$ sont exactement les mots $c_{n,j}$, où
 $0\leq j\leq d_n$. 
Tous ces mots sont des palindromes, et tout $c_{n,j}$ commence par $s_n$, sauf $c_{0,0}$ qui est vide et $c_{1,0}=a^{d_0}=c_{0,d_{0}}$.
\end{proposition}
 
\begin{remarque}{
Il ne fait aucun doute que le lien entre les palindromes sturmiens et les mots centraux est connu par les spécialistes depuis de nombreuses années. Cependant, il est difficile de trouver des références avec des déclarations précises. Par exemple, dans \cite{dlDl}, il a été démontré que chaque palindrome dans un mot sturmien est un facteur médian d’un mot central, mais cette déclaration concerne des mots et non leurs occurrences. Pour une étude récente sur les mots sturmiens bispéciaux, voir \cite{fici}.}
\end{remarque}

Ainsi, toute occurrence d'un palindrome dans un mot sturmien caractéristique ${\bf s_d}$ a son extension palindromique maximale parmi les mots $c_{n,j}$. Dans ce qui suit, nous utiliserons ce fait pour décrire les occurrences des palindromes dans ${\bf s_d}$ en termes du {\it système de numération d'Ostrowski} associé à la suite 
${\bf d}$.

Soit $q_n$ la longueur du mot standard $s_n$, c’est-à-dire
\[q_{-1} = q_0 = 1, q_{n + 1} = d_n q_n + q_{n-1} \mbox {~pour tout~} n \geq 0. \]

Dans le {\it système de numération d'Ostrowski} \cite{a_sh} associé à la suite ${\bf d} = (d_i)$, 
un entier positif $N <q_{j + 1}$ est représenté comme
\begin{equation}\label{e:numr}
N=\sum_{0\leq i \leq n} k_i q_i,
\end{equation}
où $0\leq k_i \leq d_i$ pour tout $i \geq 0$, et pour $i \geq 1$, si $k_i=d_i$, alors $k_{i-1}=0$.
Une telle représentation $N$ est unique sauf qu'on a le droit d'ajouter des zéros au début (voir Théorème 3.9.1 dans \cite{a_sh}). Dans ce chapitre, nous ne distinguerons pas les représentations qui ne diffèrent que par des zéros au début. Nous utilisons la notation $N=\overline{k_n \cdots k_1 k_0}[o]$.

Le lemme suivant est un corollaire connu des définitions d'un mot caractéristique et du système de numération d'Ostrowski.
\begin{lemme}\label{l:o}
Soit ${\bf s}$ un mot caractéristique associé à la suite directrice  ${\bf d}=(d_i)$. Considérons un nombre $N=\overline{k_n \cdots k_1 k_0}[o]$ dans le système d'Ostrowski. Alors ${\bf s}(0..N]=s_n^{k_n} s_{n-1}^{k_{n-1}}\cdots s_1^{k_1} s_0^{k_0}$.
\end{lemme}
\begin{exemple}
 {\rm
Pour la suite directrice $(1,1,1,\ldots)$, qui correspond au mot de Fibonacci,
les longueurs $q_i=|s_i|$ sont les nombres de Fibonacci, et le système de numération est le système de  Fibonacci, ou de Zeckendorf. Il décrit la décomposition gloutonne d'un nombre $N$ en somme de nombres de Fibonacci $F_n$: nous commençons ici par $F_0=1$ et $F_1=2$. Par exemple, $14=13+1=F_{5}+F_0$ et donc $14=\overline{100001}[o]$. 

Ostrowski a publié son papier sur les systèmes de numération en 1921 \cite{ostr}, mais le système de Fibonacci a été redécouvert par Zeckendorf et publié d'abord en 1952 par Lekkerkerker (qui l'attribue à Zeckendorf) \cite{lekker}, puis en 1972 par Zeckendorf lui-même \cite{zeckendorf}. 
}
\end{exemple}

Dans plusieurs articles récents (e.~g.,~\cite{efgms}), des décompositions plus générales  du même type, dites {\it légales}, ont été considérées. Une décomposition $N = \sum_{0 \leq i \leq n} k_i q_i$ est appelée {\it légale} si $ 0 \leq k_i \leq d_i$ pour $i \geq 0$, mais la deuxième restriction de la définition de la représentation d’Ostrowski n’est pas imposée. Un nombre peut admettre plusieurs représentations légales, y compris celle d'Ostrowski. Une représentation légale de $N$ est notée $N = \overline {k_n \cdots k_1 k_0}$ (sans $[o]$ à la fin, réservé à la version Ostrowski).

\begin{exemple}\label{e:14}{\rm
Dans le système de numération de Fibonacci, la représentation Ostrowski de 14 est $14=\overline{100001}[o]$, mais les représentations $14=\overline{11001}$ et $14=\overline{10111}$ sont légales aussi, car $14=F_4+F_3+F_0=8+5+1=F_4+F_2+F_1+F_0=8+3+2+1$. }
\end{exemple}

La prochaine proposition peut être trouvée par exemple dans \cite{berstel1999}.
\begin{proposition}\label{p:ddd}
Pour tout $n\geq 0$, on a $s_n^{d_n} s_{n-1}^{d_{n-1}} \cdots s_0^{d_0}=c_{n+1}$.
\end{proposition}

\begin{proposition}\label{p:lv}
Pour tous $k_0,\ldots,k_n$ tels que $k_i \leq d_i$, le mot $s_n^{k_n} s_{n-1}^{k_{n-1}} \cdots s_0^{k_0}$ est un préfixe de $c_{n+1}$, et donc de ${\bf s}$. 
\end{proposition}
 \noindent {\sc Démonstration.} 
Pour $n = 0$, l'énoncé est évident. Pour procéder par récurrence sur $n$, il suffit de prouver que si $u$ est un préfixe de $c_ {n}$, alors $s_n^{k_n}u$ est un préfixe de $c_{n + 1}$. Pour le voir, supposons d'abord que $k_n <d_n$. Dans cette situation, nous considérons $c_{n}$ comme un préfixe de $s_{n} s_{n-1}$ et voyons que $s_n^{k_n} u$ est un préfixe de $s_n^{k_n + 1} s_{n-1}$ qui est à son tour un préfixe de $s_{n + 1}$, qui est un préfixe de $c_ {n + 1}$. Supposons maintenant que $k_n = d_n$ et considérons $c_{n}$ comme un préfixe de $s_{n-1} s_ {n}$ (obtenu en effaçant les deux derniers symboles). On voit que $s_n^{k_n} u = s_n^{d_n} u$ est un préfixe de $s_n^{d_n} s_{n-1} s_n = s_{n + 1} s_n$ obtenu en effaçant au moins les deux derniers symboles, c’est-à-dire qu’il  est un préfixe de $c_{n + 1}$. \hfill $\Box$

\begin{definition}
 {
Étant données une suite directrice ${\bf d}=(d_i)$ et une suite correspondante $(s_i)$ de mots standards, $|s_i|=q_i$, on dit qu'une décomposition $N=\sum_{0\leq i \leq n} k_i q_i$ est {\it valable} si le préfixe de longueur $N$ du mot sturmien caractéristique ${\bf s_d}$ est égal à $s_n^{k_n} s_{n-1}^{k_{n-1}}\cdots s_1^{k_1} s_0^{k_0}$. Une représentation valable de $N$ est aussi notée $N=\overline{k_n \cdots k_1 k_0}$.
}
\end{definition}
Notons qu'on ne distingue pas par notation les représentations légales et valables, car, selon la proposition \ref{p:lv}, toute représentation légale est valable. Cependant, le contraire n'est pas vrai.

\begin{exemple}\label{e:1300}{\rm
Dans le système de Fibonacci, la représentation $14=\overline{1300}$ n'est pas légale,  car $k_2=3>1=d_2$. Mais elle est valable car le préfixe de longueur 14 du mot de Fibonacci est $abaababaabaaba = (abaab)(aba)^3 =s_3 s_2^3$.}
\end{exemple}

\section{Palindromes en termes d'Ostrowski}
Notre prochain objectif est le théorème suivant, améliorant légèrement le théorème 2 de \cite{frid_pal}.

\begin{theoreme}\label{t:pr}
Soit ${\bf s}$ un mot sturmien caractéristique correspondant à la suite directrice ${\bf d}=(d_n)$, et soit ${\bf s}(p_1..p_2]$ un palindrome. Alors il existe une représentation légale $p_1=\overline{x_n\cdots x_0}$ (si nécessaire, avec des zéros au début), et des nombres $m$  et $y_m$ tels que \[p_2=\overline{x_n\cdots x_{m+1} y_m \cdot (d_{m-1}-x_{m-1})\cdots (d_0-x_0)}\] est une représentation valable de $p_2$.
\end{theoreme}

\begin{exemple}
{\rm
Considérons le palindrome ${\bf f}(12..13]={\bf f}[13]=b$ dans le mot de Fibonacci
$$ {\bf f}=abaababaabaababaababaabaab\cdots.$$
Les représentations d'Ostrowski de $12=\overline{10101}[o]$ et de $13=\overline{100000}[o]$ sont toutes différentes, mais on a aussi 
$13=\overline{10110}$, donc l'énoncé du théorème est vrai avec $m=1$ et $y_m=1$.
}
\end{exemple}

Pour prouver le théorème, nous devons utiliser encore plusieurs résultats connus antérieurement, et en particulier les informations sur les occurrences de $s_n$ dans ${\bf s}$. Tout d'abord, puisque ${\bf s}$ est construit comme limite de la construction itérative \eqref{e:def}, pour tout $n$, le mot ${\bf s}$ peut être écrit comme un produit de blocs $s_n$ et $s_{n-1}$. D'après \cite{dl2002}, nous appelons cette décomposition la {\it $n$-partition} de ${\bf s}$. Nous utilisons également le lemme 3.3 du même article de Damanik et Lenz, que nous reformulons comme suit.

\begin{proposition}[\cite{dl2002}]\label{p:33}
Considérons une occurrence $s_m={\bf s}(r..r+q_m]$ de $s_m$ dans ${\bf s}$, où $m \geq 0$. Alors ${\bf s}(0..r]$ se compose de blocs complets de la  $m$-partition de ${\bf s}$, et leur suite est complètement déterminée par la lettre  ${\bf s}[r]$.
\end{proposition}

\begin{proposition}[\cite{frid_pal}]\label{p:sm}
Pour toute occurrence $s_m={\bf s}(r..r+q_m]$ d'un mot standard $s_m$ dans ${\bf s}$, où $m \geq 0$, on a
\[{\bf s}(0..r]=s_n^{k_n}\cdots s_m^{k_m}\]
de façon que
\[r=\overline{k_n\cdots k_m 00 \cdots 0}[o]\]
pour un $n \geq m$ et des $k_i\geq 0$ appropriés.
\end{proposition}

\noindent {\sc Démonstration.}[Esquisse de la preuve du théorème \ref{t:pr}] Dans \cite{frid_pal}, le théorème \ref{t:pr} a été démontré sans la condition que la décomposition de $p_1$ soit légale; cependant, la proposition \ref{p:sm} permet de l'admettre puisque la représentation nécessaire de $p_1$ peut être construite comme une concaténation de deux décompositions d'Ostrowski. En effet, considérons l'extension palindromique maximale $p_{max} = {\bf s}(p_1-l..p_2 + l]$ du palindrome $p_{ini} = 
{\bf s} (p_1..p_2]$. Comme il était discuté ci-dessus dans les propositions \ref{p:bisp} et \ref{p:central}, 
$p_{max} = c_{m, j} = s_m^j c_m$ pour certains $m$ et $j \leq d_m$; sans perte de généralité, on peut supposer que $j> 0$, puisque $c_m = c_{m, 0} = c_{m-1, d_{m-1}}$. Le palindrome initial ${\bf s}(p_1..p_2]$ est obtenu à partir de 
$p_{max}$ par effacement de $l$ symboles de gauche et $l$ symboles de droite. Soit $k$ l'entier maximal tel que 
$k q_m \leq l$, c'est-à-dire que les $k$ premières occurrences de $s_m$ à $p_{max}$ sont complètement effacées lorsque nous passons à son facteur central $p_{ini}$. Ici, $k <j$ car $c_m$ est plus court que $s_m$, et car $p_{ini}$ est situé symétriquement au centre de $p_{max}$.

Considérons maintenant le préfixe ${\bf s}(0..p_1-l + kq_m]$. Comme $k <j$, il est suivi dans ${\bf s}$ par $s_m$, et donc, selon la proposition \ref{p:sm}, la représentation d'Ostrowski de $p_1-l + kq_m$ est de la forme 
$\overline {x_n \cdots x_m 00 \cdots 0} [o] $. Maintenant ${\bf s}(0 ..p_1] = {\bf s}(0..p_1-l + kq_m] u$, où le mot $u$ est le préfixe de $s_m$ de longueur $l-kq_m <q_m$, et ainsi $u = s_{m-1}^{x_ {m-1}} \cdots s_0^{x_0}$. Si nous choisissons cette décomposition de manière gloutonne, nous obtenons la décomposition d'Ostrowski 
$l-kq_m = \overline {x_{m-1} \cdots x_0}[o]$. Ainsi, $p_1 = (p_1-l + kq_m) + (l-kq_m) = \overline {x_n \cdots x_m 00 \cdots 0}[o] + \overline {x_{m-1} \cdots x_0} [o] = \overline {x_n \cdots x_m x_m-1 \cdots x_0}$; cette dernière décomposition n'est pas toujours d'Ostrowski car nous pouvons avoir $x_m = d_m$ et $ x_ {m-1}> 0$, mais au moins elle est légale puisque pour tout $i$, $x_i \leq d_i$.

La suite de la preuve continue comme dans \cite{frid_pal}: on a $l=\overline{k x_{m-1} \cdots x_0}$, $p_2+l=\overline{x_n\cdots x_{m+1} (x_m-k+j) d_{m-1} \cdots d_0}$, et donc
\[p_2=\overline{x_n \cdots x_{m+1} \cdot (x_m-2k+j) \cdot (d_{m-1}-x_{m-1}) \cdots (d_0-x_0)}.\]
Par construction, $x_m-2k+ j\geq 0$, puisque ${\bf s}(p_1..p_2]$ est un palindrome de longueur strictement positive, il suffit de définir $x_m-2k+j=y_m$. \hfill $\Box$


\section{Conjecture sur la longueur palindromique et les mots sturmiens}
La position des palindromes par rapport au système de numération sturmien a été étudiée en relation avec la conjecture suivante, qui reste non résolue dans le cas général.

La {\it longueur palindromique} $|u|_{\rm pal}$ d'un mot fini $u$ est le plus petit nombre $Q$ de palindromes $P_1,\ldots,P_Q$ tel que $u=P_1\cdots P_Q$. 

\begin{exemple}
 {\rm
On a $|abaabb|_{\rm pal} = 3$ puisque le mot $abaabb$ ne peut pas être décomposé en concaténation de deux palindromes, mais $abaabb=(a)(baab)(b)=(aba)(a)(bb)$.
}
\end{exemple}

\begin{conjecture}\label{c:main}
Dans tout mot infini qui n'est pas ultimement périodique, la longueur palindromique des facteurs (ou même des préfixes) est non bornée.
\end{conjecture}

La conjecture a été énoncée en 2013 par Puzynina, Zamboni et l'auteur \cite{fpz} et a été démontrée dans le même document pour le cas où le mot infini évite une puissance $k$, pour un certain $k$, c'est-à-dire, ne contient pas de mot de la forme $u^k$, $u \neq \epsilon$. Il existe également une généralisation assez technique de la démonstration originale à une classe plus large de mots infinis couvrant notamment des points fixes de morphismes. Mais la preuve originale n’est pas si compliquée, et je la donne ici.

\begin{theoreme}\label{t:kpf}  Si un mot infini ${\bf w}$ évite la puissance $k$ pour un $k>0$, alors pour tout nombre naturel $P$ il existe un préfixe $u$ de ${\bf w}$ avec
$|u|_{\rm pal}>P$.
\end{theoreme}

Rappelons qu'un mot  $u=u[1]\cdots u[n]$ est appelé $t$-périodique si $u[i]=u[i+t]$ pout tout $i$ tel que $1\leq i \leq n-t$.

 La démonstration du théorème~\ref{t:kpf} utilisera les lemmes suivants.

\begin{lemme}\label{l:period}
Soit $u$ un palindrome. Alors pour tout préfixe palindrome $v$ de $u$ avec $0<|v|<|u|$, le mot $u$ est $(|u|-|v|)$-périodique.
\end{lemme}

\noindent {\sc Démonstration.}  Si $u$ et $v$ sont des palindromes, et $v$ est un préfixe de $u$, alors $v$ est aussi un suffixe de $u$ et donc $u$ est $(|u|-|v|)$-périodique.
\hfill $\Box$

\begin{lemme}\label{l:1+1/k}
Soit ${\bf w}$ un mot infini qui évite la puissance $k$. Si ${\bf w}[i_1..i_2]$ et ${\bf w}[i_1..i_3]$ sont des palindromes avec $i_3>i_2$, alors
\[\frac{|{\bf w}[i_1..i_3]|}{|{\bf w}[i_1..i_2]|}> 1+\frac{1}{k-1}.  \]
\end{lemme}
\noindent {\sc Démonstration.} Selon le lemme \ref{l:period}, le mot $w[i_1..i_3]$ est $(i_3-i_2)$-périodique; en même temps, il ne contient pas la puissance $k$, et donc $|w[i_1..i_3]|<k(i_3-i_2)$. Ainsi,
\[\frac{|{\bf w}[i_1..i_3]|}{|{\bf w}[i_1..i_2]|}=\frac{|{\bf w}[i_1..i_3]|}{|{\bf w}[i_1..i_3]|-(i_3-i_2)}>
 \frac{|{\bf w}[i_1..i_3]|}{\left (1-\frac{1}{k}\right )(|{\bf w}[i_1..i_3]|)}=1+\frac{1}{k-1}.\Box
\]
\hfill $\Box$

\begin{lemme}
Soit $N$ un nombre naturel. Alors pour tout $i\geq 0$, le nombre $n=n(i,N)$ de palindromes de la forme ${\bf w}[i..j]$  de longueur inférieure ou égale à $N$ est au plus  $2+\log_{k/(k-1)}N$.
\end{lemme}
\noindent {\sc Démonstration.} Pour tout $i\geq 0$, la longueur du plus court palindrome non vide qui commence en position $i$ est égale à $1$. Selon le lemme précédent, le prochain palindrome qui commence en position $i$ est de longueur supérieure à $\frac{k}{k-1}$, le palindrome suivant est de longueur supérieure à $(\frac{k}{k-1})^2$, et ainsi de suite. Le plus long palindrome qu'on considère est de longueur au plus $N$ et supérieure à  $(\frac{k}{k-1})^n$, ce qui implique que $n\leq \log_{k/(k-1)}N$. On voit que le nombre total $n$ de palindromes non vides qui commencent en position $i$ et sont de longueur au plus $N$ est au plus $1+\log_{k/(k-1)}N$. En ajoutant le palindrome vide nous obtenons le résultat désiré.\hfill $\Box$

\noindent {\sc Démonstration.}[Preuve du Théorème~\ref{t:kpf}]
Soit $P$ et $N$ deux entiers strictement positifs tels que \[(2+\log_{k/(k-1)}N)^P<N.\]
Par le lemme précédent, le nombre de préfixes de ${\bf w}$ de la forme $v_1v_2\ldots v_P$, où chaque $v_i$ est un palindrome, de longueur au plus $N$ est au plus $(2+\log_{k/(k-1)}N)^P$, et donc au plus $N$.
Mais ${\bf w}$ a $N$ préfixes non vides de longueur au plus $N$.
Cela signifie qu'il existe un préfixe $u$ de ${\bf w}$ de longueur au plus $N$ tel que $|u|_{\rm pal}>P$.
\hfill $\Box$

Comme nous le savons de l'article de Mignosi \cite{mignosi},  un mot sturmien évite une certaine puissance $k$ si et seulement si sa suite directrice est bornée. Dans ce cas, on peut utiliser théorème \ref{t:kpf} qui couvre par exemple le mot de Fibonacci. Cependant, le cas où la suite directrice n'est pas bornée n’a été prouvé qu’en 2018 
\cite{frid_pal}, quand nous avons réussi à utiliser le théorème \ref{t:pr} et une autre notion naturelle: la distance entre un chiffre d'une représentation et une valeur extrême pour ce chiffre.

\begin{definition}{ 
Pour une représentation valable $r=x_n\cdots x_0$ et pour tout $m\in\{0,\ldots,n\}$, notons $z_m(r)$ la distance minimale entre $x_m$ et 0 ou $x_m$ et $d_m$:
\[z_m(r)=\min(x_m,|d_m-x_m|).\]
}
\end{definition}

Avec cette notation, le théorème \ref{t:pr} a le corollaire immédiat suivant.
\begin{corollaire}\label{c:pal2repr}
Soit ${\bf w}(p_1..p_2]$ un palindrome. Alors il existe des représentations valables $r_1$ de $p_1$ et $r_2$ de $p_2$ telles que $z_i(r_1)=z_i(r_2)$ pour tout $i$ à l'exception d'une valeur $i=m$.
\end{corollaire}

Nous laissons sans preuve la proposition suivante.
\begin{proposition}[\cite{frid_pal}]\label{p:zd}
Soient $r_1$ et $r_2$ deux représentations valables du même nombre $N=\overline{r_1}=\overline{r_2}$. Alors pour tout $m$, on a $|z_m(r_1)-z_m(r_2)|\leq 3$.
\end{proposition}

\begin{exemple}
 {\rm
La différence égale à 3 est effectivement possible. Considérons une suite directrice avec $d_i=1$ pour $i=0,\ldots,3,5$ et $d_4\geq 8$. Alors
$$
\overline{140000} = \overline{1030000} = \overline{1021100} = \overline{1021011} = \overline{1020121} = \overline{1011221}.$$
On voit que le chiffre numéro 4 a diminué de 4 à 1, et donc $|z_4(140000)-z_4(1011221)|=3$.
}
\end{exemple}

\begin{theoreme}\label{t:main}
Soit ${\bf w}$ un mot sturmien caractéristique avec une suite directrice ${\bf d}$ non bornée. Alors pour tout $Q>0$, il existe un préfixe de ${\bf w}$ qui ne peut pas être décomposé en concaténation de $Q$ palindromes.
\end{theoreme}

 \noindent {\sc Démonstration.} Comme ${\bf d}$ n'est pas bornée, pour tout $Q>0$ donné, nous pouvons trouver $Q + 1$ valeurs $m_0, \cdots, m_Q$ telles que $d_{m_i} \geq 6Q + 2$ pour chaque $i$. Considérons $N = \overline {x_n \cdots x_0}$, où $x_j = 3Q + 1$ pour $n = m_i$, $i = 0, \ldots, Q$ et $x_j = 0$ sinon. Cette représentation de $N$ est légale et donc valable.

Supposons maintenant que le préfixe $w(0..N]$ puisse être représenté par une concaténation d'au plus $Q$ palindromes, c'est-à-dire qu'il existe une suite
\[0 = p_0 \leq p_1 \leq p_2 \cdots \leq p_Q = N \]
telle que pour chaque $k = 0, \ldots, Q-1$, le mot $w(p_k..p_ {k + 1}]$ est un palindrome.

En raison du corollaire \ref{c:pal2repr}, pour tout $k = 0, \ldots, Q-1$, il existe des représentations $r_{k,1}$ de $p_k$ et $r_{k + 1,2}$ de $p_{k + 1}$ telles que 
$z_i(r_{k, 1})$ et $z_i(r_{k + 1,2})$ diffèrent au maximum en un chiffre, dont nous dénotons l'indice par $i = i_k$. Notons également que 0 n'admet qu'une seule représentation valable composée de zéros (ou vide, si nous supprimons les zéros non significatifs). Par souci d'exhaustivité, nous notons la représentation $x_n \cdots x_0 $ de $N$ par $r_{Q, 1}$.

En raison de la proposition \ref{p:zd}, pour tout chiffre $m$, nous avons $|z_m(r_{k, 1}) - z_m(r_{k, 2})| \leq 3$. En même temps, la représentation 
$N = \overline{r_{Q, 1}} = \overline{x_n \cdots x_0}$ contenant $Q + 1$ chiffres $m_i$ égaux à $3Q + 1$, avec la distance $z_{m_i}(r_ {Q, 1}) = 3Q + 1$, est obtenue à partir de la représentation $0 = \overline{0 \cdots 0}$ par une série de $Q$ étapes. À chaque étape, un palindrome numéro $k$ (ici $k = 0, \ldots, Q-1$) ne peut changer que la distance $z_{i_k}(r)$ pour une valeur de $i_k$ lorsque nous passons de $r_{k , 1}$ à $r_{k + 1,2}$. Puis chaque distance $z_i(r)$ change au maximum de 3 lorsque nous passons de $r_{k + 1,2}$ à $r_{k + 1,1}$ (et de $k$ à $k + 1$).

Nous voyons qu'après $Q$ étapes de ce type, nous ne pouvons avoir que $Q$ chiffres dans la représentation dont la distance $z$ soit supérieure à $3Q$. En même temps, dans la représentation  $N = \overline{x_n \cdots x_0}$, il y en a $Q + 1$. C'est une contradiction. \hfill $\Box$

Notons que la démonstration de ce théorème ne peut pas être étendue au cas où la suite directrice ${\bf d}$ est bornée: pour ce cas-là, nous devons utiliser le théorème général \ref{t:pr} avec sa preuve «probabiliste». Notons aussi que le préfixe d'un mot caractéristique est un facteur de tout autre mot sturmien de même pente et que, grâce à un résultat élégant de Saarela \cite{saarela}, l'existence d'un préfixe de toute longueur palindromique est équivalente à l'existence d'un tel facteur. Ainsi, les deux versions de conjecture \ref{c:main} sont équivalentes et sont toutes les deux prouvées pour les mots sturmiens, même si la démonstration nécessitait deux théorèmes complètement différents.

\bibliographystyle{abbrv-fr}

\end{document}